\documentclass[manuscript,screen]{acmart}

\usepackage{algorithm}
\usepackage{algorithmicx}
\usepackage{algpseudocode}
\usepackage{hyperref}
\usepackage{listings}
\allowdisplaybreaks
\newtheorem{remark}{Remark}[section]
\usepackage{color}
\usepackage{ulem}

\AtBeginDocument{%
  \providecommand\BibTeX{{%
    \normalfont B\kern-0.5em{\scshape i\kern-0.25em b}\kern-0.8em\TeX}}}

\setcopyright{acmcopyright}
\copyrightyear{2021}
\acmYear{2021}
\acmDOI{}

\acmConference[]{}{}{}
\acmBooktitle{}
\acmPrice{}
\acmISBN{}

\begin{document}

%
\title[GCGE for Solving Large Scale Eigenvalue Problems]
{GCGE:
A Package for Solving Large Scale Eigenvalue Problems 
by Parallel Block Damping Inverse Power Method
}


\author{Yu Li}
\authornote{liyu@tjufe.edu.cn}
\email{liyu@tjufe.edu.cn}
\orcid{0000-0002-2246-8005}
\affiliation{%
	\\
	\institution{Coordinated Innovation Center for Computable Modeling in Management Science, Tianjin University of Finance and Economics, Tianjin, 300222, China}
}

\author{ZiJing Wang}
\authornote{zjwang@lsec.cc.ac.cn}
\email{zjwang@lsec.cc.ac.cn}
\author{Hehu Xie}
\authornote{hhxie@lsec.cc.ac.cn}
\email{hhxie@lsec.cc.ac.cn}
\affiliation{%
	\\
	\institution{LSEC, ICMSEC, Academy of Mathematics and Systems Science, Chinese Academy of Sciences, Beijing, 100190, China}
	\\
	\institution{School of Mathematical Sciences, University of Chinese Academy of Sciences, Beijing, 100049, China}
}


\renewcommand{\shortauthors}{Li, Wang and Xie}

\begin{abstract}
We propose an eigensolver and the corresponding package, GCGE, for solving large scale eigenvalue problems.  This method is the combination of damping idea, subspace projection method and inverse power method with dynamic shifts.  To reduce the dimensions of projection subspaces, a moving mechanism is developed when the number of desired eigenpairs is large.  The numerical methods, implementing techniques and the structure of the package are presented.  Plenty of numerical results are provided to demonstrate the efficiency, stability and scalability of the concerned eigensolver and the package GCGE for computing many eigenpairs of large symmetric matrices arising from applications.
\end{abstract}



\maketitle

\section{Introduction}

A fundamental and challenging task in modern science and engineering is to solve large scale
eigenvalue problems. Although high-dimensional eigenvalue problems are ubiquitous in physical sciences, data and
imaging sciences, and machine learning, there is no so many classes of eigensolvers
as that of linear solvers. Compared with linear equations, there are less efficient numerical
methods for solving large scale eigenvalue problems, which poses significant challenges for scientific computing
\cite{Bai2000Templates}.
In particular, the eigenvalue problems from
complicated systems bring strong demand for eigensolvers with good
efficiency, stability and scalability 
at the same time
\cite{Fan2014Parallel,Fan2015Parallel,Yu2018Parallel}.

The Krylov subspace methods such as Arnoldi and Lanczos methods 
are always used to design the eigensolvers
\cite{Saad1992Numerical}.
In order to use explicitly and implicitly restarted techniques for generalized eigenvalue problems,
it is necessary to solve the included linear equations exactly to produce upper Hessenberg matrices.
But this requirement is always very difficult for large scale sparse matrices with poor conditions.
Based on this consideration, 
LOBPCG is designed based on some types of
iteration processes which do not need to solve the included linear equations exactly
\cite{Knyazev2003Efficient,Knyazev2006Toward,
Hetmaniuk2006Basis,Knyazev2007Block,Duersch2018Robust}.
This property makes LOBPCG be reasonable candidate 
for solving large scale eigenvalue problems on parallel computers.
But the subspace generating method and orthogonalization way 
lead to the unstability of LOBPCG algorithm \cite{Li2020parallel,Zhang2020generalized}.

The appearance of high performance computers
brings new issues for computing plenty of eigenpairs of large scale matrices, 
which has not so good efficiency and scalability as solving large scale linear equations. 
Solving eigenvalue problems on high performance computers 
needs new considerations about the stability and scalability 
of orthogonalization for plenty of vectors, 
efficiency and memory costing for computing Rayleigh-Ritz problems.
The aim of this paper is to develop a method and 
the corresponding package 
for solving symmetric eigenvalue problems. 
This method
is the combination of damping idea, subspace projection method 
and inverse power method with dynamic shifts.
The package GCGE (Generalized Conjugate Gradient Eigensolver) 
is written by C language and constructed with the way of matrix-free and vector-free.
In order to improve the efficiency, stability and scalability, 
we also introduce new efficient implementing
techniques for orthogonalization and computing Rayleigh-Ritz problems. 
A recursive orthogonalization method with SVD (Singular Value Decomposition) is proposed
in order to improve parallel efficiency.
In addition,
we also provide a moving mechanism to reduce 
the dimensions of projection subspaces when solving Rayleigh-Ritz problems.
The source code can be downloaded from
GitHub with the address \url{https://github.com/Materials-Of-Numerical-Algebra/GCGE}.




The rest of the paper is organized as follows. 
In Section \ref{Section_GCGE}, we present the concerned algorithm for eigenvalue problems.
The implementing techniques are designed in Section \ref{Section_Optimization}.
In Section \ref{Section_Numerical_Example}, plenty of numerical
tests are provided to demonstrate the efficiency, stability and scalability of the proposed algorithm 
and the associated package.
Concluding remarks are given in the last section.

\section{GCG algorithm}\label{Section_GCGE}
For simplicity, in this paper, we are concerned with the following generalized algebraic eigenvalue problem:
Find eigenvalue $\lambda\in\mathbb R$ and eigenvector $x\in\mathbb R^N$ such that
\begin{eqnarray}\label{Eigenvalue_Problem}
Ax = \lambda B x,
\end{eqnarray}
where $A$ is $N\times N$ real symmetric matrix
and $B$ is $N\times N$ real symmetric positive definite (SPD) matrix.

The generalized conjugate gradient (GCG) algorithm is a type of subspace projection method, which uses the block damping inverse power idea to
generate triple  blocks $[X,P,W]$, where $X$ saves the current eigenvector approximation,
$P$ saves the information from previous iteration step,
and $W$ saves vectors from $X$ by the inverse power iteration with some CG steps.
We refer to this method as generalized conjugate gradient algorithm since the structure of triple blocks $[X,P,W]$
is similar to that of conjugate gradient method.
Assuming that it is desired to compute the smallest {\tt numEigen} eigenpairs, 
the corresponding GCG algorithm is defined by Algorithm \ref{GCG_Algorithm},
where {\tt numEigen} stands for the number of desired eigenpairs.

\begin{algorithm}[!hb]
	\caption{GCG algorithm} \label{GCG_Algorithm}
\begin{algorithmic}[0]
\State 1. Choose {\tt numEigen} vectors to build the block $X$ and two null blocks $P=[\ ]$, $W=[\ ]$.
\State 2. Define $V=[X,P,W]$ and do orthogonalization to $V$ in the sense of inner product deduced by the matrix $B$.
\State 3. Solve the Rayleigh-Ritz problem 
$(V^{\top}AV)\hat{x}= \hat{x}\Lambda_x$ 
to obtain $\Lambda_x$ and $\hat{x}$, then
get new approximate eigenvectors $X^{\tt  new} = V\hat{x}$.
\State 4. Check the convergence of eigenpair approximations $(\Lambda_x, X^{\tt  new})$.
If the smallest {\tt numEigen} {eigenpairs converge}, the iteration will stop.
\State 5. Otherwise, compute $P = X^{\tt  new}-X(X^{\top}BX^{\tt new})$
and update $X=X^{\tt new}$.
\State 6. {Generate} $W$ by solving linear equations
$(A-\theta B)W = BX(\Lambda_x-\theta I)$ by some CG steps with the initial guess $X$,
where the shift $\theta$ is selected dynamically.
\State 7. Then go to {\bf STEP 2}.
\end{algorithmic}
\end{algorithm}

The main difference of Algorithm \ref{GCG_Algorithm} from LOBPCG is the way to generate $W$ and orthogonalization to $V$.
The GCG algorithm uses the inverse power method with dynamic shifts to generate $W$.
Meanwhile, the full orthogonalization to $V$ is implemented 
in order to guarantee the numerical stability.
In addition, 
a new type of recursive orthogonalization method with SVD is designed 
in the next section.

In {\bf Step 3} of Algorithm \ref{GCG_Algorithm}, solving Rayleigh-Ritz problem is a sequential process
which can not be accelerated by using normal parallel computing.
Furthermore, it is well known that the computing time is superlinearly dependent on the number of
desired eigenpairs \cite{Saad1992Numerical}.
Then in order to accelerate this part, reducing the dimensions of Rayleigh-Ritz problems is a reasonable way.
We will compute the desired eigenpairs in batches when the number of desired eigenpairs is large.
In each iteration step, 
the dimensions of $P$ and $W$ are set to be 
${\tt numEigen}/5$ or ${\tt numEigen}/10$.
Moreover,
a moving mechanism is presented
for computing large number of desired eigenpairs.
These two strategies can further not only reduce
the time proportion of the sequential process
for solving Rayleigh-Ritz problems
but also reduce the amount of memory
required by {\bf STEP 3}.
In addition, the Rayleigh-Ritz problem is 
distributed to multi computing processes
and each process only computes a small part of desired eigenpairs.
In other words, the Rayleigh-Ritz problem is solved in parallel.
More details of implementing techniques will be introduced
in Sections \ref{sec:reduce_computation}
and \ref{sec:moving_mechanism}.

In {\bf STEP 6} of  Algorithm \ref{GCG_Algorithm}, 
though the matrix $A-\theta B$ may not be SPD,
the CG iteration method is adopted for solving the included linear equations
due to the warm start $X$ and the shift $\theta$.
Please see Section \ref{sec:reduce_computation} for more details.
Furthermore, it is suggested to use the algebraic multigrid method as the preconditioner 
for {\bf STEP 6} of Algorithm \ref{GCG_Algorithm} with the shift $\theta=0.0$, 
when the concerned matrices are sparse and come
from the discretization of partial differential operators
by finite element, finite difference or finite volume, etc.

\section{Implementing techniques}\label{Section_Optimization}

In this section, we introduce implementing techniques to improve efficiency, scalability and stability for the concerned eigensolver in this paper.
Based on the discussion in the previous sections, 
we focus on the methods for doing the orthogonalization and
computing Rayleigh-Ritz problems.
A recursive orthogonalization method with SVD and a moving mechanism are presented.
In addition,
the package GCGE is introduced,
which is written by C language and constructed with the way of matrix-free and vector-free.


\subsection{Improvements for orthogonalization}\label{GCGE:subsec:orthogonal}

This subsection is devoted to introducing the orthogonalization methods
which have been supported by GCGE.
So far, we have provided
modified block orthogonalization method and recursive orthogonalization method with SVD.
The criterion for choosing the orthogonalization methods should be based on the number of desired eigenpairs
and the scales of the concerned matrices. 
The aim is to keep the balance among efficiency, stability and scalability.



The modified Gram-Schmidt method \cite{Stewart2008Block} 
is designed to improve the stability of classical orthogonalization method.
The modified block orthogonalization method 
is the block version of modified Gram-Schmidt method,
which can be defined by Algorithm \ref{Modified_Block_Orth}.
They have the same accuracy and stability,
but the modified block orthogonalization method has better efficiency and scalability.

Let us consider the orthogonalization for $X\in\mathbb{R}^{N\times m}$ 
and assume $m=b\ell$ in Algorithm \ref{Modified_Block_Orth}.
We divide $X$ into $\ell$ blocks, i.e.,
$X=[X_1,X_2,\cdots,X_\ell]$, where
$X_i\in\mathbb{R}^{N\times b}$, $i=1,\cdots,\ell$.
The orthogonalization process is to make $X$ be orthogonal
to $X_0$ and do orthogonalization for $X$ itself,
where $X_0\in\mathbb{R}^{N\times m_0}$ has already been orthogonalized, i.e., $X_0^{\top}BX_0=I$.

Firstly, in order to maintain the numerical stability, the process of deflating components in $X_0$ from $X$ is
repeated until the norm of $X_0^{\top}BX$ is small enough.
Secondly,	
the columns of $X$ in blocks of $\ell$ columns are orthogonalized through the modified Gram-Schmidt method.
For each $k=1,\cdots,\ell$ in Algorithm \ref{Modified_Block_Orth},
when $X_k$ is linear dependent, 
the rearmost vectors of $X$ are copied to the corresponding location.
In addition, Algorithm \ref{Modified_Block_Orth} needs 
$b+1$ global communications in each ${\tt for}$ iteration.
In other words, the total number of global communications is
\begin{equation*}
	(b+1)(\ell-1)+b=m+m/b-1.
\end{equation*}

In fact, in modified block orthogonalization method, we deflate the components in previous orthogonalized vectors
successively for all unorthogonalized vectors in each iteration step.
This means Algorithm \ref{Modified_Block_Orth}
uses block treatment for the unorthogonalized
vectors to improve efficiency and scalability without loss of stability.
As default, $b$ is set to be $\min(m/4,200)$.

\begin{algorithm}[!htb]
\caption{Modified block orthogonalization} \label{Modified_Block_Orth}
\begin{algorithmic}[1]
\Repeat
	\State Compute
	$X = X - X_0 (X_0^{\top}(BX))$;
\Until {the norm of $X_0^{\top}(BX)$ is small enough;}
   \For { $k = 1 : \ell$}
	 \State Orthogonalize $X_k$ by modified Gram-Schmidt method;
	 \If {$k==\ell$} \State break; \EndIf
\Repeat
       \State Compute
				$\begin{bmatrix} R_{k+1} \\ \vdots \\ R_{\ell} \end{bmatrix}
				= [{X}_{k+1},\cdots,{X}_{\ell}]^{\top} (BX_{k})$;
       \State Compute
		 $[X_{k+1},\cdots,X_{\ell}] = [{X}_{k+1},\cdots,{X}_{\ell}]
				 -{X}_{k} \begin{bmatrix}R_{k+1} \\ \vdots \\ R_{\ell} \end{bmatrix}^{\top}$;
				 \Until {the norm of $R_{k+1},\cdots,R_{\ell}$ are small enough;}
   \EndFor
\end{algorithmic}
\end{algorithm}


In order to improve efficiency and scalability further,
we design a type of recursive orthogonalization method with SVD and the corresponding scheme is defined by Algorithm \ref{RecusiveOrthSVD}.
The aim here is to take full use of level-3 BLAS operations.
We also find the paper \cite{Yokozawa2006Efficient} has discussed the similar orthogonalization method without SVD.
The contribution here is to combine the recursive orthogonalization method and SVD to improve the scalability.

\begin{algorithm}[!htb]
	\caption{${\tt RecusiveOrthSVD}($X$,{\tt s},{\tt e})$} \label{RecusiveOrthSVD}
	\begin{algorithmic}[1]
	\State Compute ${\tt length} = {\tt e}-{\tt s}+1$;
		\If{${\tt length} \leq c$}
		\Repeat
		\State Compute $M=X({\tt :~,s:e})^{\top}BX({\tt :~,s:e})$;
		\State Compute SVD of $M=Q\Lambda Q^{\top}$;
		\State Compute $X({\tt :~,s:e})=X({\tt :~,s:e})Q\Lambda^{-1/2}$;
		\Until {the norm of $\Lambda-I$ is small enough;}
		\Else
	\State ${\tt s1}={\tt s}$; ${\tt e1}={\tt s}+{\tt length}/2-1$;
\State ${\tt s2} = {\tt e1}+1$; ${\tt e2} = {\tt e}$;
	\State Call ${\tt RecusiveOrthSVD}(X,{\tt s1},{\tt e1})$;
\Repeat
	\State Compute $R=(BX({\tt :~,s1:e1}))^{\top}X({\tt :~,s2:e2})$;
	\State Compute 
	$X({\tt :~,s2:e2})
	=X({\tt :~,s2:e2})-X({\tt :~,s1:e1})R$;
				 \Until {the norm of $R$ are small enough;}
	\State Call ${\tt RecusiveOrthSVD}(X,{\tt s2},{\tt e2})$;
		\EndIf
	\end{algorithmic}
\end{algorithm}

Let us consider $X\in\mathbb{R}^{N\times m}$ and $m=2^\eta$
in Algorithm \ref{RecusiveOrthSVD}.
The orthogonalization of $X$ is completed 
by calling ${\tt RecusiveOrthSVD}$ recursively.
We use $X({\tt :~,s:e})$ to stand for the {\tt s}-th column to {\tt e}-th column of $X$.
When ${\tt length}\leq c$, SVD is applied to 
computing $X$, 
where $c$ is set to be $\min(m,16)$ as default.
In order to maintain the numerical stability, 
computing $X$ with SVD
is repeated until 
the matrix $\Lambda$ is close to the identity matrix.
In general, the above condition is satisfied
after two or three iterations.
If $M$ has eigenvalues close to zero, i.e., the set of vectors is linearly dependent,
the subsequent vectors will be copied to the corresponding location.

{If $c=16$ and we compute $X$ with SVD three times when ${\tt length}\leq c$, }
the total number of global communications is
\begin{equation*}
	2^0+2^1+2^2+\cdots+2^{\eta-5}+3\times2^{\eta-4} =
	\frac{1}{4}m-1,
\end{equation*}
which is much less than the total number of global communications of Algorithm \ref{Modified_Block_Orth}.



The recursive orthogonalization method with SVD is recommended and it is the default choice in our package for the orthogonalization to long vectors.
In fact, Algorithms \ref{Modified_Block_Orth} and \ref{RecusiveOrthSVD} can both reach the required accuracy for all numerical examples in this paper.
In the case of solving generalized eigenvalue problems, $B$-orthogonalization should be considered.
Algorithm \ref{RecusiveOrthSVD} is more efficient than Algorithm \ref{Modified_Block_Orth} {in most cases},
which will be shown in Section \ref{sec:nev_large}.


\subsection{Computation reduction for Algorithm \ref{GCG_Algorithm}}\label{Subsection_Algorithm_Optimization}
\label{sec:reduce_computation}

In this subsection, let us continue considering the whole computation procedure for Algorithm \ref{GCG_Algorithm}.
The aim here is to design
efficient ways to compute the Rayleigh-Ritz problem in {\bf STEP 3} which include
\begin{itemize}
	\item Orthogonalizing to $V=[X,P,W]$;
	\item Computing the small scale matrix $\bar{A} = V^{\top}AV$;
	\item Solving the standard eigenvalue problem $\bar{A} \hat{x}=\hat{x}\Lambda_x$.
\end{itemize}

Except for the moving mechanism shown 
in Section \ref{sec:moving_mechanism} and the inverse power method with 
dynamic shifts for solving $W$,
the techniques here are almost the same as that
in \cite{Li2020parallel, Zhang2020generalized}.
But for easier understanding and
completeness, we also introduce them here
using more concise expressions.
In conclusion, the following main optimization techniques are implemented:
\begin{itemize}
	\item[(1)] The converged eigenpairs do not participate the subsequent iteration;
	\item[(2)] The sizes of $P$ and $W$ are set to be {\tt blockSize}, which is equal to {\tt numEigen/5} as default;
	\item[(3)] The shift is selected dynamically when solving $W$;
	\item[(4)] The large scale orthogonalization to $V$ is transformed into the small scale orthogonalization to $P$
		and a large scale orthogonalization to $W$;
	\item[(5)] The submatrix of $\bar{A}$ corresponding to $X$ can be obtained by $\Lambda_x$;
	\item[(6)] The submatrix of $\bar{A}$ corresponding to $P$ can be computed by multiplication of small scale dense matrices;
	\item[(7)] The Rayleigh-Ritz problem $\bar{A} \hat{x}=\hat{x}\Lambda_x$ is solved in parallel;
  \item[(8)] The moving mechanism is presented to reduce the dimension of $\bar{A}$ further.
\end{itemize}

According to {\bf STEP 2}  of Algorithm \ref{GCG_Algorithm}, 
we decompose $X$ into three parts
\begin{equation*}
	X =\begin{bmatrix} X_c, & X_n, & X_{\widetilde{n}}\end{bmatrix},
\end{equation*}
where $X_c$ denotes the converged eigenvectors and $[X_n, X_{\widetilde{n}}]$  denotes the unconverged ones.
The number of vectors in  $X_n$ is {\tt blockSize}.
Based on the structure of $X$,
the block version has the following structure
\begin{equation*}
	V=\begin{bmatrix}
		X_c, & X_n, & X_{\widetilde{n}}, & P, & W
\end{bmatrix}
\end{equation*}
with $V^{\top}BV=I$.
And the eigenpairs $\Lambda_x$ and $\hat{x}$ can be decomposed into the following form
\begin{equation}
	\label{equ:def_lambda_x}
\Lambda_x =
\begin{bmatrix}
\Lambda_c & O & O\\
O &\Lambda_n & O\\
O & O & \Lambda_{\widetilde{n}}
\end{bmatrix},\
\hat{x} = \begin{bmatrix}
	\hat{x}_c, & \hat{x}_n, & \hat{x}_{\widetilde{n}}
\end{bmatrix},
\end{equation}
where $\Lambda_x$ is the diagonal matrix.

Then in {\bf STEP 3} of Algorithm \ref{GCG_Algorithm}, the small scale eigenvalue problem 
\begin{equation*}
\bar{A}\hat{x} = \hat{x}\Lambda_x
\end{equation*}
has the following form
\begin{equation}\label{equ:step3_old}
\bar{A}
\begin{bmatrix}
I &O & O\\
O &\hat{x}_{nn} & \hat{x}_{n\widetilde{n}}\\
O &\hat{x}_{\widetilde{n}n} & \hat{x}_{\widetilde{n}\widetilde{n}}\\
O &\hat{x}_{pn}& \hat{x}_{p\widetilde{n}}\\
O &\hat{x}_{wn}& \hat{x}_{w\widetilde{n}}\\
\end{bmatrix}
=
\begin{bmatrix}
I &O & O\\
O &\hat{x}_{nn} & \hat{x}_{n\widetilde{n}}\\
O &\hat{x}_{\widetilde{n}n} & \hat{x}_{\widetilde{n}\widetilde{n}}\\
O &\hat{x}_{pn}& \hat{x}_{p\widetilde{n}}\\
O &\hat{x}_{wn}& \hat{x}_{w\widetilde{n}}\\
\end{bmatrix}
\begin{bmatrix}
\Lambda_c & O & O\\
O &\Lambda_n & O\\
O & O & \Lambda_{\widetilde{n}}
\end{bmatrix},
\end{equation}
where $\bar{A} = V^{\top}AV$, $\hat{x}^{\top}\hat{x}=I$ and
$\hat{x}_c$, $\hat{x}_n$, $\hat{x}_{\widetilde{n}}$ have following structures
\begin{equation}\label{Structure_XPW}
\hat{x}_c =
\begin{bmatrix}
I \\ O\\ O\\ O\\ O
\end{bmatrix},\
\hat{x}_n =
\begin{bmatrix}
O\\ \hat{x}_{nn}\\ \hat{x}_{\widetilde{n}n}\\ \hat{x}_{pn}\\ \hat{x}_{wn}
\end{bmatrix},\
\hat{x}_{\widetilde{n}} =
\begin{bmatrix}
O \\
\hat{x}_{n\widetilde{n}}\\
\hat{x}_{\widetilde{n}\widetilde{n}}\\
\hat{x}_{p\widetilde{n}}\\
\hat{x}_{w\widetilde{n}}
\end{bmatrix}.
\end{equation}
In addition, Ritz vectors is updated as
\begin{equation*}
	X^{\tt  new}=V\hat{x}.	
\end{equation*}


In {\bf STEP 4} of Algorithm \ref{GCG_Algorithm},
the convergence of the eigenpairs $(\Lambda_x, X^{\tt  new})$ is checked.
Due to (\ref{equ:def_lambda_x}),
we set 
\begin{align*}
	&\hat{x}_n=\begin{bmatrix} \hat{x}_{n_1,} & \hat{x}_{n_2} \end{bmatrix},\
	\hat{x}_{\widetilde{n}}= \begin{bmatrix} \hat{x}_{\widetilde{n}_1}, & \hat{x}_{\widetilde{n}_2} \end{bmatrix},\\
&\Lambda_n=\begin{bmatrix}
	\Lambda_{n_1} & O \\
	O &\Lambda_{n_2}
\end{bmatrix},\
\Lambda_{\widetilde{n}}=\begin{bmatrix}
	\Lambda_{\widetilde{n}_1} & O \\
	O &\Lambda_{\widetilde{n}_2}
\end{bmatrix},
\end{align*}
and 
the diagonal of $\Lambda_{n_1}$ inclues the new converged eigenvalues.
Then
all $\ell$ convergened eigenvectors are in
\begin{equation*}
	X^{\tt{new}}_c= V\begin{bmatrix} \hat{x}_c, & \hat{x}_{n_1} \end{bmatrix},\
\end{equation*}
and 
the unconverged ones are in
\begin{equation*}
	X^{\tt{new}}_n= V\begin{bmatrix} \hat{x}_{n_2}, & \hat{x}_{\widetilde{n}_1} \end{bmatrix} \mbox{ and }\
	X^{\tt{new}}_{\widetilde{n}}= V\hat{x}_{\widetilde{n}_2}.
\end{equation*}
If $\ell$ is equal to {\tt numEigen}, the iteration will stop.
Otherwise, $0\leq \ell<{\tt numEigen}$ and the iteration will continue.
Here, 
the length of $X^{\tt{new}}_n$ is
\begin{equation*}
	{\tt blockSize} = \min({\tt numEigen/5}, {\tt numEigen}-\ell).
\end{equation*}

In {\bf STEP 5} of Algorithm \ref{GCG_Algorithm}, in order to produce $P$ for the next GCG iteration,
from the definition of $\hat{x}_n$ in (\ref{Structure_XPW})
and the orthonormality of $V$, i.e., $V^{\top}BV=I$, we first set
\begin{equation*}
	\widetilde{P} =V\hat{x}_n - X_n(X_n^{\top}BV\hat{x}_n)
	=V\widetilde{p},
\end{equation*}
where
\begin{equation}\label{CP_tilde}
\widetilde p=
\begin{bmatrix}
O \\ O \\ \hat{x}_{\widetilde{n}{n}}\\ \hat{x}_{pn}\\ \hat{x}_{wn}
\end{bmatrix}.
\end{equation}

In order to compute $P^{\tt new}$ to satisfy $(X^{\tt new})^{\top}BP^{\tt new}=O$,
we come to do the orthogonalization for small scale vectors in
$[\hat{x}, \widetilde p]$
according to the $L^2$ inner product.
Since vectors in $\hat{x}$ are already orthonormal,
the orthogonalization only needs to be done
for $\widetilde p$ against $\hat{x}$ to get a new vectors $\hat{p}$.
Thus let
$[\hat{x}, \hat{p}]$
denote the orthogonalized block, i.e.,
\begin{equation}
	\label{Equality_Cpx}
	\begin{bmatrix}\hat{x}, & \hat{p}\end{bmatrix}^{\top}
	\begin{bmatrix}\hat{x}, & \hat{p}\end{bmatrix} = I.
\end{equation}
Then,
\begin{equation*}
	P^{\tt{new}} = V\hat{p}.
\end{equation*}
Moreover, it is easy to check that
\begin{equation*}
	(X^{\tt{new}})^{\top}BP^{\tt{new}} = \hat{x}^{\top}V^{\top}BV\hat{p} = O
\end{equation*}
and
\begin{equation*}
	(P^{\tt{new}})^{\top}BP^{\tt{new}} = \hat{p}^{\top}V^{\top}BV\hat{p} = I.
\end{equation*}

In {\bf STEP 6} of Algorithm \ref{GCG_Algorithm},
$\widetilde{W}$ is obtained by some CG iterations
for the linear equations
\begin{equation}
	\label{equ:compW}
(A-\theta B)\widetilde{W}
= B X^{\tt{new}}_{n}
\begin{bmatrix}\Lambda_{n_2}-\theta I&O\\
O&\Lambda_{\widetilde{n}_1}-\theta I\end{bmatrix}
\end{equation}
with the initial guess $X^{\tt{new}}_{n}$, where
the shift $\theta$ is set to be the largest converged eigenvalue in the convergence process.
It is noted that the shift is not fixed and the matrix $A-\theta B$ may not be SPD,
but the initial guess $X^{\tt{new}}_{n}$ is perpendicular 
to the eigenvectors of $A-\theta B$ corresponding to all negtive eigenvalues,
i.e.,
\begin{equation*}
	(X^{\tt{new}}_n)^{\top}(A-\theta B)X^{\tt{new}}_c = O,
\end{equation*}
since $X^{\tt{new}}_c$ reaches the convergence criterion.
In other words, $A-\theta B$ is SPD in the orthogonal complement space of $\mbox{span}(X^{\tt{new}}_c)$.
Then the CG iteration method can be adopted for solving the included linear equations.
Due to the shift $\theta$,
the multiplication of matrix and vector of each CG iteration
takes more time,
but the convergence of GCG algorithm is accelerated.
In addition, there is no need to solve linear equations 
(\ref{equ:compW})
with high accuracy,
and only $10$-$30$ CG iterations are enough
during each GCG iteration.
In Remark	\ref{rem:shift}, an example is presented to explain
why the convergence of GCG algorithm with dynamic shifts 
is accelerated after one CG iteration.
In Section \ref{sec:cg_max_niter},
we give some numerical results to show the performance of GCGE with 
dynamic shifts and 
the convergence procedure under different number of CG iterations.
In order to produce  $W^{\tt{new}}$ for the next GCG iteration, 
using Algorithm \ref{RecusiveOrthSVD},
we need to
do the orthognalization to $\widetilde{W}$ according to $[X^{\tt{new}}, P^{\tt{new}}]$, i.e.,
\begin{equation*}
\begin{bmatrix}
	X^{\tt{new}}, & P^{\tt{new}}
\end{bmatrix}^{\top}
B W^{\tt{new}} = O,\
(W^{\tt{new}})^{\top}B W^{\tt{new}} = I.
\end{equation*}

\begin{remark}
	\label{rem:shift}
	{We give an example to present the accelerating convergence of GCG algorithm with dynamic shifts after one CG iteration.}
	Assuming that the first eigenpair $(\lambda_1,v_1)$ has been found for the standard eigenvalue problem
	\begin{equation*}
		Ax=\lambda x,
	\end{equation*}
	we have the approximate eigenvector $x_0 = a_2v_2+a_3v_3$
	of the second eigenvector,
	where 
	\begin{equation*}
	Av_2=\lambda_2v_2,\	Av_3=\lambda_3v_3, \mbox{ and }\ 0<\lambda_1<\lambda_2\leq\lambda_3.
	\end{equation*}
For the linear equations
\begin{equation*}
	(A-\theta I)w=(\tilde{\lambda}-\theta)x_0 \mbox{ and }\ 0\leq\theta<\lambda_2,
\end{equation*}
we can obtain the new approximate eigenvector
	\begin{equation*}
		x_1 = \frac{a_2^2+a_3^2}{a_3^2(\lambda_2-\theta)+a_2^2(\lambda_3-\theta)}\Big((\lambda_3-\theta)a_2v_2+(\lambda_2-\theta)a_3v_3\Big),
	\end{equation*}
after the first CG iteration 
with 
the initial guess $x_0$,
where $\tilde{\lambda}=x_0^{\top}Ax_0/x_0^{\top}x_0$. 
{
It is noted that the convergence rate is
	\begin{equation*}
		\frac{\lambda_2-\theta}{\lambda_3-\theta},
	\end{equation*}
which is less than the case of $\theta=0$.}
\end{remark}

Backing to {\bf STEP 2} of Algorithm \ref{GCG_Algorithm}, we denote
\begin{equation*}
	V^{\tt{new}}=
\begin{bmatrix}
	X^{\tt{new}}, & P^{\tt{new}}, & W^{\tt{new}}
\end{bmatrix}=
\begin{bmatrix}
	V\hat{x}, & V\hat{p}, & W^{\tt{new}}
\end{bmatrix}.
\end{equation*}
During solving the Rayleigh-Ritz problem,
we need to assemble the small scale matrices
$(V^{\tt{new}})^{\top}AV^{\tt{new}}$ and $(V^{\tt{new}})^{\top}BV^{\tt{new}}$.
Since the orthogonalization to the vectors in $V^{\tt{new}}$
has been done by the inner product deduced by the matrix $B$,
$(V^{\tt{new}})^{\top}BV^{\tt{new}}$ is an identity matrix.
Then we only need to compute the matrix $\bar{A}_{\tt  new}$,
which is equal to
\begin{equation}\label{Equality_A_new}
\begin{bmatrix}
(X^{\tt{new}})^{\top}AX^{\tt{new}} & (X^{\tt{new}})^{\top}AP^{\tt{new}} & (X^{\tt{new}})^{\top}AW^{\tt{new}}\\
	(P^{\tt{new}})^{\top}AX^{\tt{new}} & (P^{\tt{new}})^{\top}AP^{\tt{new}} & (P^{\tt{new}})^{\top}AW^{\tt{new}}\\
	(W^{\tt{new}})^{\top}AX^{\tt{new}} & (W^{\tt{new}})^{\top}AP^{\tt{new}} & (W^{\tt{new}})^{\top}AW^{\tt{new}}
 \end{bmatrix}.
\end{equation}

From (\ref{equ:step3_old}), the submatrix $(X^{\tt{new}})^{\top}AX^{\tt{new}}$ does not need to
be computed explicitly since {it satisfies the following formula}
\begin{equation}\label{Equality_xx}
(X^{\tt{new}})^{\top}AX^{\tt{new}} =
\hat{x}^{\top}V^{\top}AV\hat{x} = \hat{x}^{\top}\bar{A} \hat{x}
= \Lambda_x.
\end{equation}
Based on the basis in $V$ and (\ref{Equality_Cpx}),
we have 
\begin{equation}\label{Equality_pp}
(P^{\tt{new}})^{\top}AP^{\tt{new}} =
\hat{p}^{\top}V^{\top}AV\hat{p} = \hat{p}^{\top}\bar{A} \hat{p}
\end{equation}
and
\begin{equation}\label{Equality_xp}
(P^{\tt{new}})^{\top}AX^{\tt{new}} =
\hat{p}^{\top}V^{\top}AV\hat{x} = \hat{p}^{\top}\bar{A} \hat{x}
= \hat{p}^{\top} \hat{x} \Lambda_x = O.
\end{equation}
Thus from (\ref{Equality_A_new}), (\ref{Equality_xx}), (\ref{Equality_pp}) and (\ref{Equality_xp}),
we know the matrix $\bar{A}^{\tt{new}}$ has the following structure
\begin{equation}\label{equ:struct_bar_A_new}
\begin{bmatrix}
\Lambda_0 & O & O& O\\
O&	\Lambda_1 & O & \alpha_1\\
O&	O   & \alpha_0 & \alpha_2 \\
O&	\alpha_1^{\top} & \alpha_2^{\top} & \alpha_3 \\
 \end{bmatrix},
\end{equation}
where
\begin{align*}
&\Lambda_0=
\begin{bmatrix}
	\Lambda_c & O \\ O & \Lambda_{n_1}
\end{bmatrix},\
\Lambda_1=
\begin{bmatrix}
	\Lambda_{n_2} & O & O\\
O & \Lambda_{\widetilde{n}_1} & O\\
O & O & \Lambda_{\widetilde{n}_2}
\end{bmatrix},\\
&\alpha_0=\hat{p}^{\top}\bar{A}\hat{p},\
\alpha_1=
	\begin{bmatrix}
		X^{\tt new}_n, & X^{\tt{new}}_{\widetilde{n}}
	\end{bmatrix}^{\top} AW^{\tt{new}},\\
&\alpha_2=
(P^{\tt{new}})^{\top} AW^{\tt{new}},\
\alpha_3=
(W^{\tt{new}})^{\top} AW^{\tt{new}}.
\end{align*}
It is noted that
since $X^{\tt new}_c$ reaches the convergence criterion,
we assume the equation
\begin{equation*}
AX^{\tt new}_c = BX^{\tt new}_c\Lambda^{\tt new}_c
\end{equation*}
is satisfied.
Then
\begin{equation*}
(W^{\tt new})^{\top}AX^{\tt new}_c = (W^{\tt new})^{\top} BX_c^{\tt new}\Lambda^{\tt new}_c = O
\end{equation*}
is satisfied approximately since $(W^{\tt new})^{\top}BX^{\tt new}_c=O$.

After assembling matrix $\bar A^{\tt new}$, the next task is to solve
the new small scale eigenvalue problem:
\begin{equation}
\label{Small_Eigenvalue_Problem_new}
	\bar{A}^{\tt{new}} \hat{x}^{\tt{new}} =
\hat{x}^{\tt{new}} \Lambda_x^{\tt{new}},
\end{equation}
in {\bf STEP 3}.
Due to the converged eigenvectors $X^{\tt new}_c$ in $V^{\tt new}$,
there are already $\ell$ converged eigenvectors of $\bar{A}^{\tt{new}}$
and they all have the form 
\begin{equation*}
(0,...,0,1,0,...0)^{\top}	
\end{equation*}
($1$ stays in the position of associated converged eigenvalue).
We only need to compute the unconverged eigenpairs corresponding to $[X^{\tt new}_n, X^{\tt new}_{\widetilde n}]$
for the eigenvalue problem (\ref{Small_Eigenvalue_Problem_new}).
The subroutine {\tt dsyevx} from LAPACK \cite{Anderson1999LAPACK} is called to compute
the only $(\ell+1)$-th to ${\tt numEigen}$-th eigenvalues and
their associated eigenvectors.

In order to reduce time consuming of this part,
this task is distributed to multi computing processes and
each process only computes a small part of desired eigenpairs.
After all processes finish their tasks, the subroutine ${\tt  MPI\_Allgatherv}$ is adopted
to gather all eigenpairs from all processes and deliver them to all.
This way leads to an obvious time reduction for computing the desired eigenpairs of
(\ref{Small_Eigenvalue_Problem_new}).
Since more processes lead to more communicating time, we choose the number of
used processes for solving (\ref{Small_Eigenvalue_Problem_new})
such that each process computes at least $10$ eigenpairs.


\begin{remark}
In order to accelerate the convergence,
the size of $X$, {\tt sizeX}, is always chosen to be greater than {\tt numEigen}, 
which is set to be 
the minimum of ${\tt numEigen}+3\times{\tt blockSize}$ and the dimension of $A$,
as default.
\end{remark}

\begin{remark}
	Since the converged eigenpairs $(\Lambda_c,X_c)$ do not participate in the subsequent iterations, 
in real implementation,
$\bar{A}$ is computed as follows
\begin{equation*}
\begin{bmatrix}
	X_n, & X_{\widetilde{n}}, & P, & W
\end{bmatrix}^{\top}A
\begin{bmatrix}
	X_n, & X_{\widetilde{n}}, & P, & W
\end{bmatrix},
\end{equation*}
and the corresponding eigenpairs have the forms
\begin{equation*}
\begin{bmatrix} \Lambda_n & O \\ O& \Lambda_{\widetilde{n}} \end{bmatrix},\
\begin{bmatrix}
 \hat{x}_{nn} &\hat{x}_{n\widetilde{n}} \\
 \hat{x}_{\widetilde{n}n} &\hat{x}_{\widetilde{n}\widetilde{n}}\\
 \hat{x}_{pn} &\hat{x}_{p\widetilde{n}}\\
 \hat{x}_{wn} &\hat{x}_{w\widetilde{n}}
 \end{bmatrix}.
\end{equation*}
	In other words, the internal locking (deflation) is implemented to 
prevent computing over again the eigenpairs which have been found.
\end{remark}

\subsection{The moving mechanism}
\label{sec:moving_mechanism}

In Algorithm  \ref{GCG_Algorithm},
the small scale eigenvalue problem (\ref{Small_Eigenvalue_Problem_new})
needs to be solved,
in which the dimension of the dense matrix $\bar{A}$ is
$${\tt sizeX}+2\times{\tt blockSize},$$
where the size of $X$, {\tt sizeX}, is equal to ${\tt numEigen}+3\times{\tt blockSize}$.  
When {\tt numEigen} is large, e.g., $5000$, with ${\tt blockSize}=200$,
{\tt dsyevx} should be called to solve $5000$ eigenpairs for a dense matrix of $6000$-dimension.
In this case,  the time of {\bf STEP 3} of Algorithm \ref{GCG_Algorithm} is always dominated.

\begin{figure}[!htb]
\centering
\includegraphics[width=8cm,height=3cm]{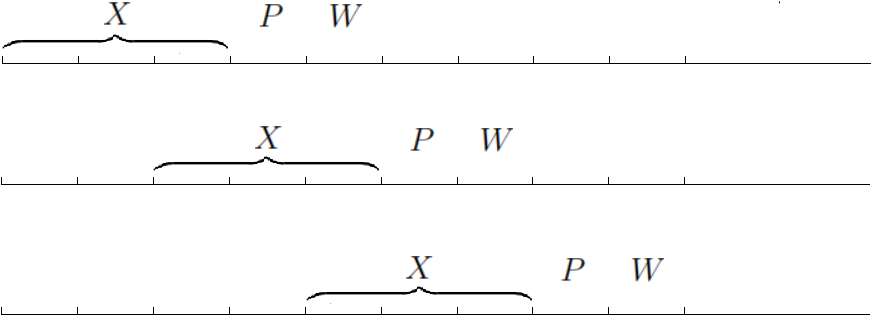}
\caption{Moving $[X,P,W]$, when $2\times{\tt blockSize}$ eigenpairs converged.}
\Description{Moving $[X,P,W]$, when $2\times{\tt blockSize}$ eigenpairs converged.}
\label{fig:xpw}
\end{figure}

In order to improve efficiency further for the above case,
we present a moving mechanism. 
Firstly, 
the maximum project dimension is set to be ${\tt maxProjDim}=5\times{\tt blockSize}$
in moving procedure,
i.e., the size of $X$ is set to be $3\times{\tt blockSize}$ 
and the sizes of $P$ and $W$ are both ${\tt blockSize}$.
Secondly, when $2\times{\tt blockSize}$ eigenpairs converged,
all the eigenpairs 
of $\bar{A}$ will be solved, i.e,
$\bar{A}$ is decomposed into 
\begin{equation*}
	\bar{A}=\begin{bmatrix}\hat{x},&\hat{p},&\hat{w}\end{bmatrix}
	\Lambda_{xpw}\begin{bmatrix}\hat{x},&\hat{p},&\hat{w}\end{bmatrix}^{-1},
\end{equation*}
where
\begin{equation*}
	\bar{A}=V^{\top}AV,\
	V = \begin{bmatrix}X,&P,&W
	\end{bmatrix}.
\end{equation*}
In addition, 
the new $X$ is equal to $V[\hat{x},\hat{p},\hat{w}]$, 
and $\Lambda_{xpw}$ can be used to construct the new $\bar{A}$
in the next {\bf STEP 3}.
In other words, $P$ and $W$ have been integrated into $X$.   
Then, the new $P$ and $W$ will be computed and stored behind the new $X$.
When there are new converged $2\times{\tt blockSize}$ eigenpairs again,
$P$ and $W$ will be integrated into $X$ again, and so on.
The above process is shown in Figure \ref{fig:xpw}.
It is noted that 
the dimension of the dense matrix $\bar{A}$ is ${\tt maxProjDim}=5\times{\tt blockSize}$ at most
in the small scale eigenvalue problem (\ref{Small_Eigenvalue_Problem_new}).

\begin{figure}[!htb]
\centering
\includegraphics[scale=0.45]{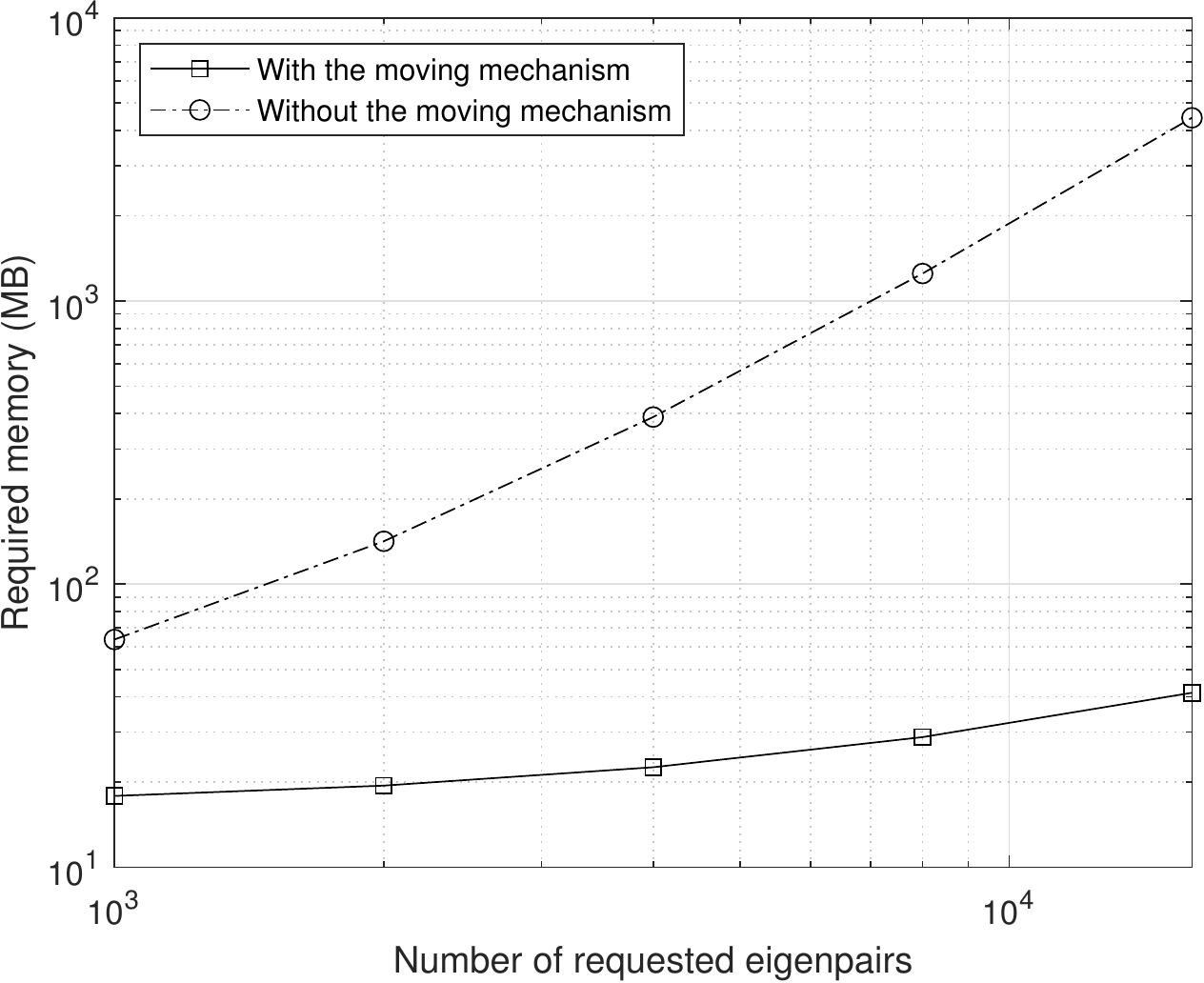}
\caption{Requested memory in each process}
\Description{}
\label{fig:moving_mem}
\end{figure}

Moreover, the moving mechanism can greatly reduce memory requirements,
which allows more eigenpairs to be computed.
Specifically speaking, 
the double array, of which the size is
\begin{equation}
	\begin{split}
		&({\tt sizeX}+2\times{\tt blockSize})+2\times({\tt maxProjDim})^2\\
		&+10\times({\tt maxProjDim})+{\tt sizeX}\times{\tt blockSize},
	\end{split}
	\label{equ:moving_mem}
\end{equation}
is required to be stored in each process.
The first two terms denote the sizes of the two arrays 
which are used to store the eigenpairs and the dense matrix
in the small scale eigenvalue problem (\ref{Small_Eigenvalue_Problem_new}).
The third term is the size of workspace for {\tt dsyevx}.
The last term is the size of the array which is used in {\bf STEP 5. }
In Figure \ref{fig:moving_mem}, the required memory computed by (\ref{equ:moving_mem}) is shown
with and without the moving mechanism.

\subsection{Matrix-free and vector-free operations}\label{sec:mat_vec_free}


Based on Algorithm \ref{GCG_Algorithm} and its implementing techniques presented in above sections, 
we develop the package GCGE, which 
is written by C language and constructed
with the way of matrix-free and vector-free.
So far, {the package has included}
the eigensolvers for the matrices which are stored in
dense format, compressed row/column sparse format or are supported in MATLAB,
Hypre \cite{Falgout2006design}, PETSc \cite{Balay1997Efficient}, PHG \cite{Zhang2009parallel}
and SLEPc \cite{Hernandez2005SLEPc}.
Table \ref{matrix-vector-free} presents the currently supported matrix-vector structure.
It is noted that there is no need 
to copy the built-in matrices and the vectors 
from these softwares/libraries to the GCGE package.

\begin{table}[!htb]
\centering
\caption{Supported matrix-vector structures}\label{matrix-vector-free}
\begin{tabular}{c|cc}
	\toprule
       & matrix structure name    & vector structure name\\     
	\midrule
MATLAB & sparse distributed matrix & full stored matrix \\
Hypre  & hypre\_ParCSRMatrix & hypre\_ParVector \\
PETSc  & Mat                 & Vec              \\
PHG    & MAT                 & VEC              \\
SLEPc  & Mat                 & BV               \\
	\bottomrule
\end{tabular}
\end{table}

A user can also build his own eigensolver by 
providing the matrix, vector structures and their operations.
The following six matrix-vector operations
should be provided by the user:
\begin{itemize}
	\item[(1)] {\tt VecCreateByMat}
	\item[(2)] {\tt VecDestroy}
	\item[(3)] {\tt VecLocalInnerProd}
	\item[(4)] {\tt VecSetRandomValue}
	\item[(5)] {\tt VecAxpby}
	\item[(6)] {\tt MatDotVec}
\end{itemize}
They realize creating and destroying vector according to matrix,
computing local inner product of vectors $x$ and $y$,
setting random values for vector $x$,
computing vector $y = \alpha x+\beta y$,
computing vector $y = Ax$,
respectively.
{\tt VecInnerProd}, i.e.,
computing inner product of vectors $x$ and $y$,
has been provided
through calling
{\tt VecLocalInnerProd}
and {\tt MPI\_Allreduce}.


The default matrix-multi-vector operations are invoked based on the above matrix-vector operations and the additional two operations:
{\tt GetVecFromMultiVec} and {\tt RestoreVecForMultiVec},
which are getting/restoring one vector from/to multi vectors.
For higher efficiency, it is strongly recommended that
users should provide the following six matrix-multi-vector operations:
\begin{itemize}
	\item[(1)] {\tt MultiVecCreateByMat}
	\item[(2)] {\tt MultiVecDestroy}
	\item[(3)] {\tt MultiVecLocalInnerProd}
	\item[(4)] {\tt MultiVecSetRandomValue}
	\item[(5)] {\tt MultiVecAxpby}
	\item[(6)] {\tt MatDotMultiVec}
\end{itemize}
In addition, if user-defined multi-vector is stored in dense format,
BLAS library can be used to implement
(1)-(5) operators easily,
which has been provided in the GCGE package.
In other words,
only one operator, i.e.,
computing the multiplication of matrix and multi-vector
needs to be provided by users.


In order to improve the parallel efficiency of
computing inner products of multi-vector $X$ and $Y$,
i.e., the operation {\tt MultiVecInnerProd},
a new MPI data type with the corresponding reduced operation
has been created by
\begin{verbatim}
            MPI_Type_vector, MPI_Op_create.
\end{verbatim}
The variable {\tt MPI\_IN\_PLACE} is used as the value of {\tt sendbuf} in
{\tt MPI\_Allreduce}
at all processes.

Although SLEPc \cite{Hernandez2005SLEPc} provides an inner product operation for BV structure,
we still recommend using our own multi-vector inner product operation.
Let us give an example to illustrate the reason.
For instance, we need to compute the inner products
\begin{equation*}
	[x_i,\cdots,x_j]^{\top}[y_p,\cdots,y_q]	
\end{equation*}
and the results are stored in
the following submatrix
\begin{equation}\label{Dense_SubMatrix}
	\begin{bmatrix}
	  c_{ip} & \cdots & c_{iq}\\
	  \vdots &   & \vdots\\
		c_{jp} & \cdots & c_{jq}
 \end{bmatrix}.
\end{equation}
Always, the vectors $[x_i,\cdots,x_j]$ and $[y_p,\cdots,y_q]$ come from the multi-vector
\begin{align*}
	X&=[x_1,\cdots,x_i,\cdots,x_j,\cdots,x_n],\\
	Y&=[y_1,\cdots,y_p,\cdots,y_q,\cdots,y_m].
\end{align*}
and the dense matrix (\ref{Dense_SubMatrix}) is one submatrix of the following matrix
\begin{equation*}
	\begin{bmatrix}
	 * & * & \cdots & * &  *\\
	 * & c_{ip} & \cdots & c_{iq} &  *\\
	 * & \vdots &  & \vdots &  *\\
	 * & c_{jp} & \cdots & c_{jq} &  *\\
	 * & * & \cdots & * &  *
 \end{bmatrix}_{s\times t}
\end{equation*}
which is stored by column. Thus, it can be noted that the above mentioned submatrix (\ref{Dense_SubMatrix})
is not stored continuously.

The result of the SLEPc's inner product operation, {\tt BVDot}, 
must be stored in a sequential dense matrix with dimensions $n\times m$ at least.
In other words, 
regardless of the values of $i$, $j$, $p$ and $q$,
in each process, the additional memory space is required,
of which the size is $n\times m$.
In general, $n$ and $m$ are set to be ${\tt sizeX}+2\times{\tt blockSize}$ in the GCG algorithm, 
while $s$ and $t$ are much less than $n$ and $m$, respectively.

In the GCGE package, the operation {\tt MultiVecInnerProd} is implemented 
as follows:
\begin{itemize}
	\item[(1)]
Through {\tt MultiVecLocalInnerProd},
local inner products are calculated and stored in
the above mentioned submatrix for each process;
	\item[(2)]
A new {\tt MPI\_Datatype} named {\tt SUBMAT} is created by
\begin{verbatim}
int MPI_Type_vector(
    int count, int length, int stride,
    MPI_Datatype oldtype, MPI_Datatype *newtype)
\end{verbatim}
with
\begin{verbatim}
    count=q-p+1, length=j-i+1, stride=s;
\end{verbatim}
	\item[(3)]
Through {\tt MPI\_Op\_create},
the operation of sum of {\tt SUBMAT} is created,
which is named as {\tt SUM\_SUBMAT};
	\item[(4)] Then 
\begin{verbatim}
int MPI_Allreduce(
    void *sendbuf, void *recvbuf, int count,
    MPI_Datatype datatype, MPI_Op op,
    MPI_Comm comm)
\end{verbatim}
is called with
\begin{verbatim}
        sendbuf=MPI_IN_PLACE, count=1,
        datatype=SUBMAT, op=SUM_SUBMAT
\end{verbatim}
to gather values from all processes
and distribute the results back to all processes.
\end{itemize}
Obviously, no extra workspace is needed here.
The memory requirements are reduced for each process.

\section{Numerical results}\label{Section_Numerical_Example}

The numerical experiments in this section are carried out on LSSC-IV
in the State Key Laboratory of Scientific and Engineering Computing,
Chinese Academy of Sciences.
Each computing node has two 18-core Intel Xeon Gold 6140 processors at 2.3 GHz and 192 GB memory.
For more information, please check
\url{http://lsec.cc.ac.cn/chinese/lsec/LSSC-IVintroduction.pdf}.
We use {\tt numProc} to denote the number of processes in numerical experiments.

In this section, the GCG algorithm defined by Algorithm \ref{GCG_Algorithm} 
and the implementing techniques in Section \ref{Section_Optimization} are
investigated for thirteen standard eigenvalue problems and one generalized eigenvalue problem.
The first thirteen matrices
are available in
Suite Sparse Matrix Collection\footnote {\url{https://sparse.tamu.edu}}, 
which have clustered eigenvalues and many negative eigenvalues.
The first matrix named Andrews is provided by Stuart Andrews at Brown University, which has seemingly random sparsity pattern.
The second to the thirteenth matrices are generated by
the pseudo-potential algorithm for real-space electronic structure calculations 
\cite{Kronik2006PARSEC, Natan2008Real, Saad2010Numerical}.
The FEM matrices $A$ and $B$ come from the finite element discretization for the following Laplace eigenvalue problem:
Find $(\lambda,u)\in \mathbb R\times H_0^1(\Omega)$ such that
\begin{eqnarray}\label{Laplace_Eigenvalue_Problem}
\left\{
\begin{array}{rcl}
-\Delta u &=&\lambda u,\ \ \ \ {\rm in}\ \Omega,\\
u&=&0,\ \ \ \ \ \ {\rm on}\ \partial\Omega,
\end{array}
\right.
\end{eqnarray}
where $\Omega=(0,1)\times (0,1)\times (0,1)$.
The discretization of the eigenvalue problem (\ref{Laplace_Eigenvalue_Problem})
by the conforming cubic finite element (P3 element) with 3,145,728 elements
leads to the stiffness matrix $A$ and the mass matrix $B$.
The concerned matrices are listed in Table \ref{tab:matrix},
where 
the density is defined by
\begin{equation*}
\frac{\mbox{the number of non-zero entries}}{\mbox{dimension}\times\mbox{dimension}}.
\end{equation*}
The proposed GCG algorithm given by Algorithm \ref{GCG_Algorithm} based on BV structure from SLEPc is adopted
to solve eigenpairs of the concerned matrices in Table \ref{tab:matrix}.

\begin{table}[!htb]
\centering
\caption{Testing matrices}\label{tab:matrix}
\begin{tabular}{clrrr}
	\toprule
ID&Matrix       & Dimension  & Non-zero Entries & Density     \\     
	\midrule                   
1&Andrews       & 60,000     & 760,154     &  2.11{\tt e}-4   \\
2&{CO}          & 221,119    & 7,666,057   &  1.57{\tt e}-4   \\
3&Ga10As10H30   & 113,081	   & 6,115,633   &  4.78{\tt e}-4   \\
4&{Ga19As19H42} & 133,123	   & 8,884,839   &  5.01{\tt e}-4   \\ 
5&Ga3As3H12     & 61,349     & 5,970,947   &  1.59{\tt e}-3   \\
6&Ga41As41H72   & 268,096    & 18,488,476  &  2.57{\tt e}-4   \\
7&{Ge87H76}     & 112,985    & 7,892,195   &  6.18{\tt e}-4   \\ 
8&{Ge99H100}    & 112,985    & 8,451,395   &  6.62{\tt e}-4   \\ 
9&{Si34H36}     & 97,569     & 5,156,379   &  5.42{\tt e}-4   \\
10&{Si41Ge41H72}& 185,639    & 15,011,265  &  4.36{\tt e}-4   \\
11&Si5H12       & 19,896     & 738,598     &  1.87{\tt e}-3   \\
12&{Si87H76}    & 240,369    & 10,661,631  &  1.85{\tt e}-4   \\
13&SiO2         & 155,331    & 11,283,503  &  4.68{\tt e}-4   \\
14&FEM matrices $A$ and $B$ & 14,045,759 & 671,028,055 &  3.40{\tt e}-6 \\
	\bottomrule
\end{tabular}
\end{table}

The convergence criterion is set to be
\begin{equation*}
	\|Ax-\lambda x\|_2/\|x\|_2 < {\tt  tol}
\end{equation*} 
for the first thirteen matrices and 
\begin{equation*}
	\|Ax-\lambda Bx\|_2/(\lambda\|B^{1/2}x\|_2) < {\tt tol}
\end{equation*} 
for FEM matrices,
where the tolerance, {\tt tol}, is set to be $10^{-8}$ as default.
Moreover, we set 
${\tt blockSize}={\tt numEigen}/10$ 
for the first thirteen matrices and 
${\tt blockSize}={\tt numEigen}/5$ 
for FEM matrices.

In order to confirm  the efficiency, stability and scalability of GCGE,
we investigate the numerical comparison between GCGE and LOBPCG.
We will find that GCGE has better efficiency, stability than LOBPCG 
and they have almost the same scalability.
In addition, Krylov-Schur method is also compared
in Sections \ref{sec:tol_small} and \ref{sec:nev_large}.

\subsection{{About dynamic shifts and the number of CG iterations}}
\label{sec:cg_max_niter}

	In this subsection, we give some numerical results to show the performance
of GCGE with dynamic shifts 
and the convergence procedure under different number of CG iterations.

\begin{figure}[!htb]
\centering
\includegraphics[scale=0.45]{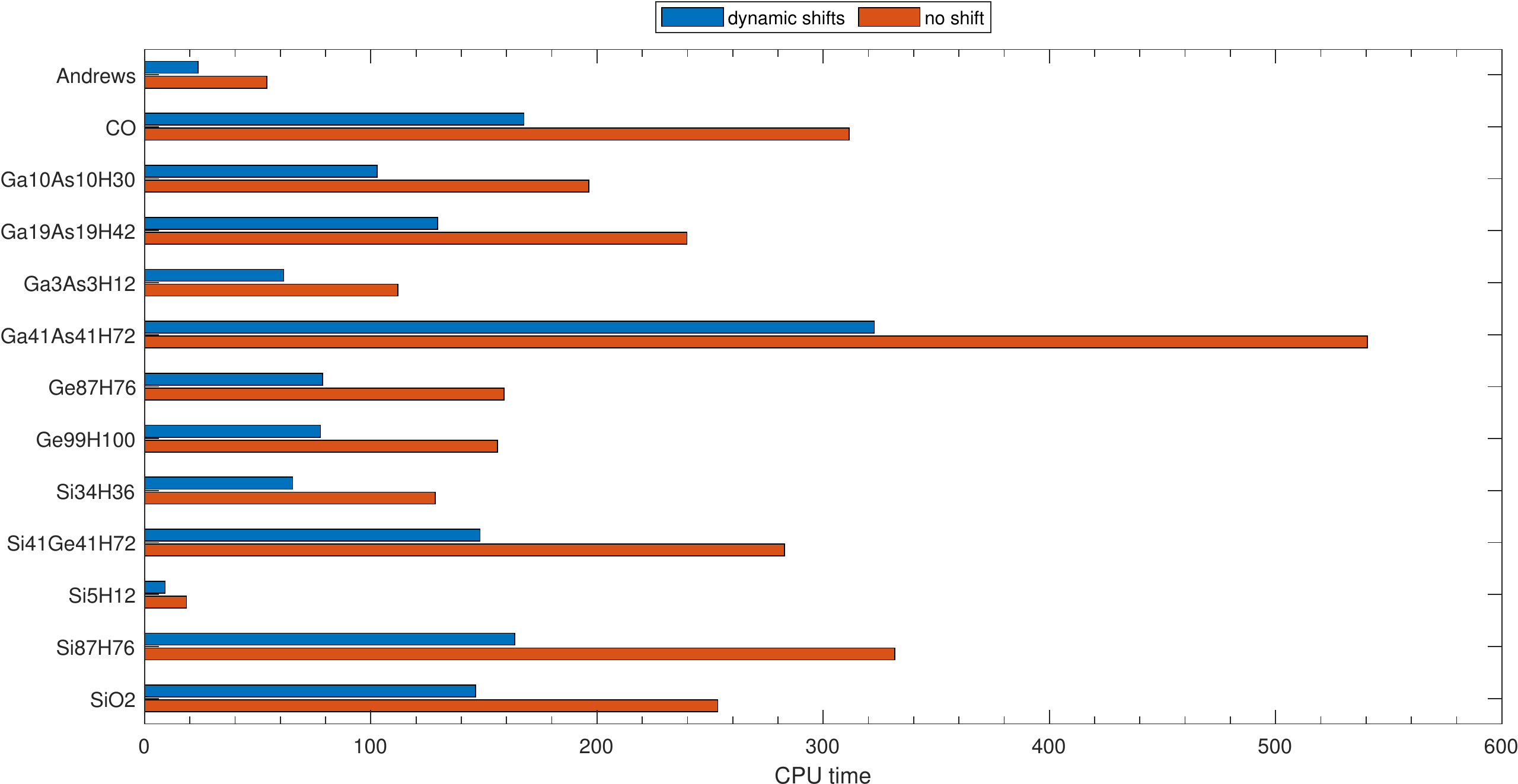}
\caption{${\tt tol}=10^{-8}$, ${\tt numEigen}=800$, and ${\tt numProc}=36$}
\label{fig:shifted_1e-8}
\end{figure}

In {\bf STEP 6} of Algorithm \ref{GCG_Algorithm},
the linear equations (\ref{equ:compW})
are solved by some CG iterations.
Due to the shift $\theta$,
the multiplication of matrix and vector of each CG iteration
takes more time,
but the convergence of GCG algorithm is accelerated.
For the standard eigenvalue problems, i.e., $B=I$, 
because the additional computation is only the linear operations on vectors, 
each GCG iteration with dynamic shifts takes a little more time than the case of no shift.
As shown in Figure \ref{fig:shifted_1e-8},
the performance of GCGE with dynamic shifts is greatly improved.
In addition, the total number of GCG iterations is presented in Table \ref{tab:avg_time_shift}.


\begin{table}[!htb]
\centering
\caption{The total number of GCG iterations}
\begin{tabular}{c|rrrrr}
	\toprule
ID&Matrix       & Dynamic Shifts & No Shift& Ratio \\     
	\midrule                   
1&Andrews       & 102 & 281 & 36.29\% \\
2&{CO}          & 97  & 195 & 49.74\% \\
3&Ga10As10H30   & 105 & 213 & 49.29\% \\
4&{Ga19As19H42} & 110 & 216 & 50.92\% \\ 
5&Ga3As3H12     & 81  & 165 & 49.09\% \\
6&Ga41As41H72   & 133 & 236 & 56.35\% \\
7&{Ge87H76}     & 78  & 212 & 36.79\% \\ 
8&{Ge99H100}    & 77  & 206 & 37.37\% \\ 
9&{Si34H36}     & 79  & 207 & 38.16\% \\
10&{Si41Ge41H72}& 87  & 208 & 41.82\% \\
11&Si5H12       & 86  & 201 & 42.78\% \\
12&{Si87H76}    & 89  & 232 & 38.36\% \\
13&SiO2         & 90  & 164 & 54.87\% \\
\bottomrule                            
\end{tabular}                          
	\label{tab:avg_time_shift}               
\end{table}

For the generalized eigenvalue problems,
there is no significant improvement for 
the overall performance of GCGE with dynamic shifts 
by the additional computation of the multiplication of matrix $B$ and vectors.
When the matrix $A$ can be modified, 
we recommend users to generate $A-\theta B$ explicitly and 
do CG steps for $A-\theta B$ directly. 
In this event, GCGE with dynamic shifts will perform better
for the generalized eigenvalue problem 
and the results for ${\tt numEigen}=800$ and {\tt numEigen = 5000}
are shown in Tables \ref{tab:shift_FEM} and \ref{tab:nev5000}
respectively.

\begin{table}[!htb]
\centering
\caption{FEM matrices with ${\tt numEigen}=800$, ${\tt tol}=10^{-12}$ and ${\tt numProc}=576$}
\begin{tabular}{l|rr}
	\toprule
                & The Total Number  & CPU Time (in seconds) \\     
                & of GCG Iterations &                       \\     
	\midrule                   
	Dynamic Shifts  & 83 & 1669.19 \\
	No Shift        & 88 & 1777.87 \\
	Ratio           & 94.31\% & 93.88\% \\
\bottomrule                            
\end{tabular}                          
	\label{tab:shift_FEM}               
\end{table}

In addition, the GCG algorithm do not need to solve linear equations exactly in {\bf STEP 6}.
In the rest of this subsection, 
the total time of GCGE and the average time per each GCG iteration are presented 
under different number of CG iterations.
Because the first thirteen matrices have similar density, 
we choose SiO2 with ${\tt numProc}=36$ 
and FEM matrices with ${\tt numProc}=576$
as test objects.

\begin{figure}[!htb]
\centering
\includegraphics[scale=0.45]{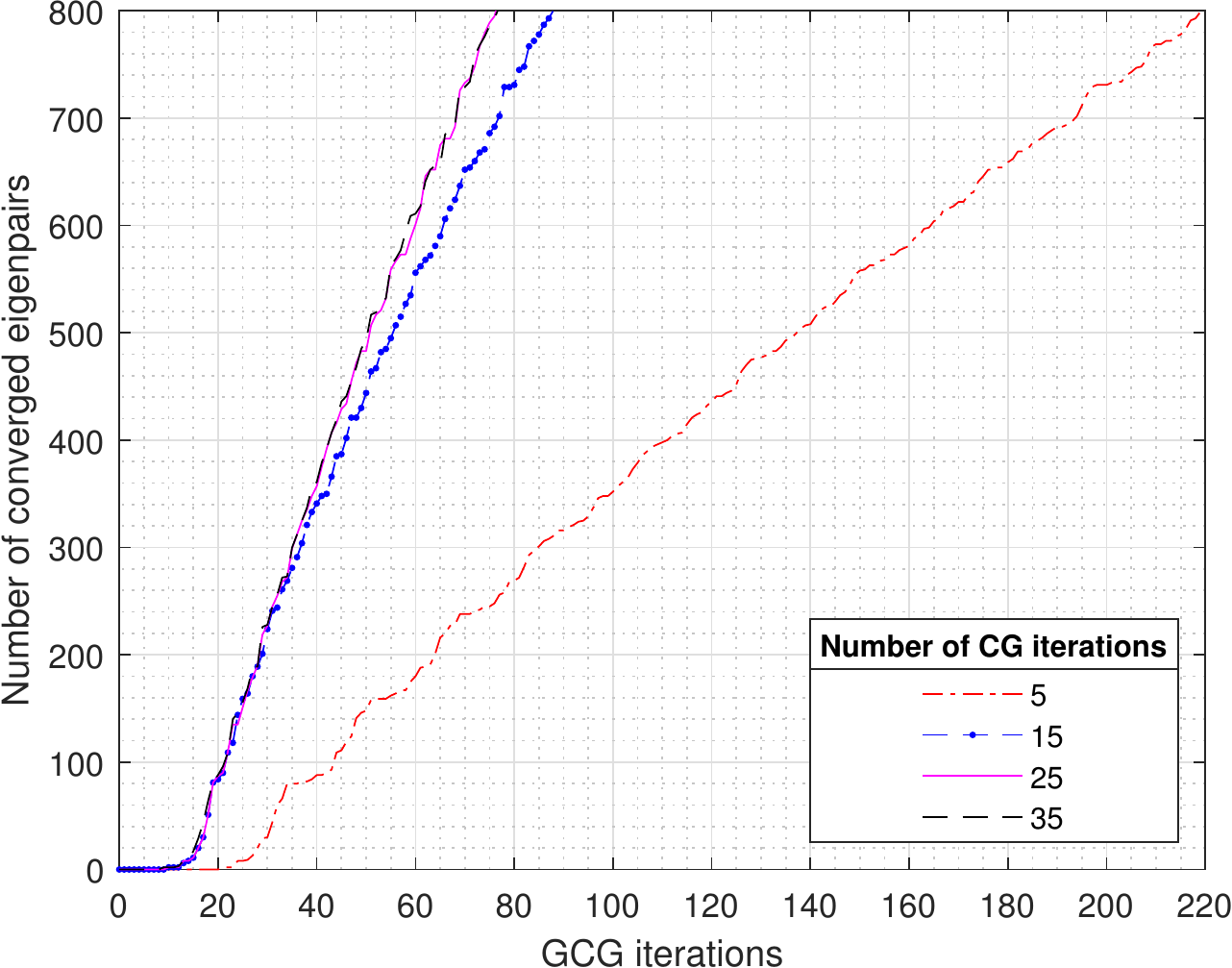}
\caption{Convergence procedure for SiO2 with ${\tt tol}=10^{-8}$}
\label{fig:cg_nconv_SiO2}
\end{figure}

\begin{figure}[!htb]
\centering
\includegraphics[scale=0.45]{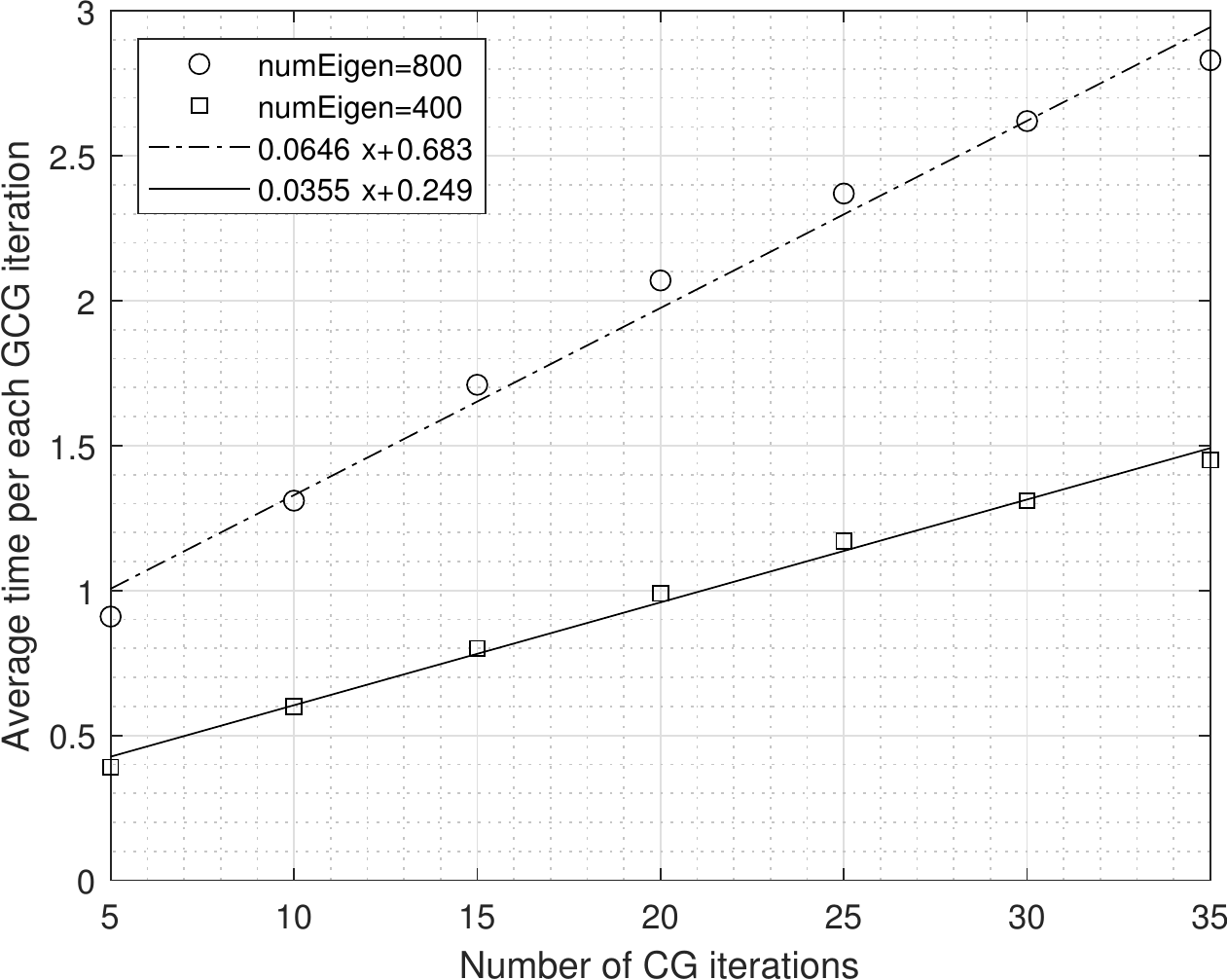}
\caption{Average time for SiO2 with ${\tt tol}=10^{-8}$}
\label{fig:cg_avg_time_SiO2}
\end{figure}

\begin{figure}[!htb]
\centering
\includegraphics[scale=0.45]{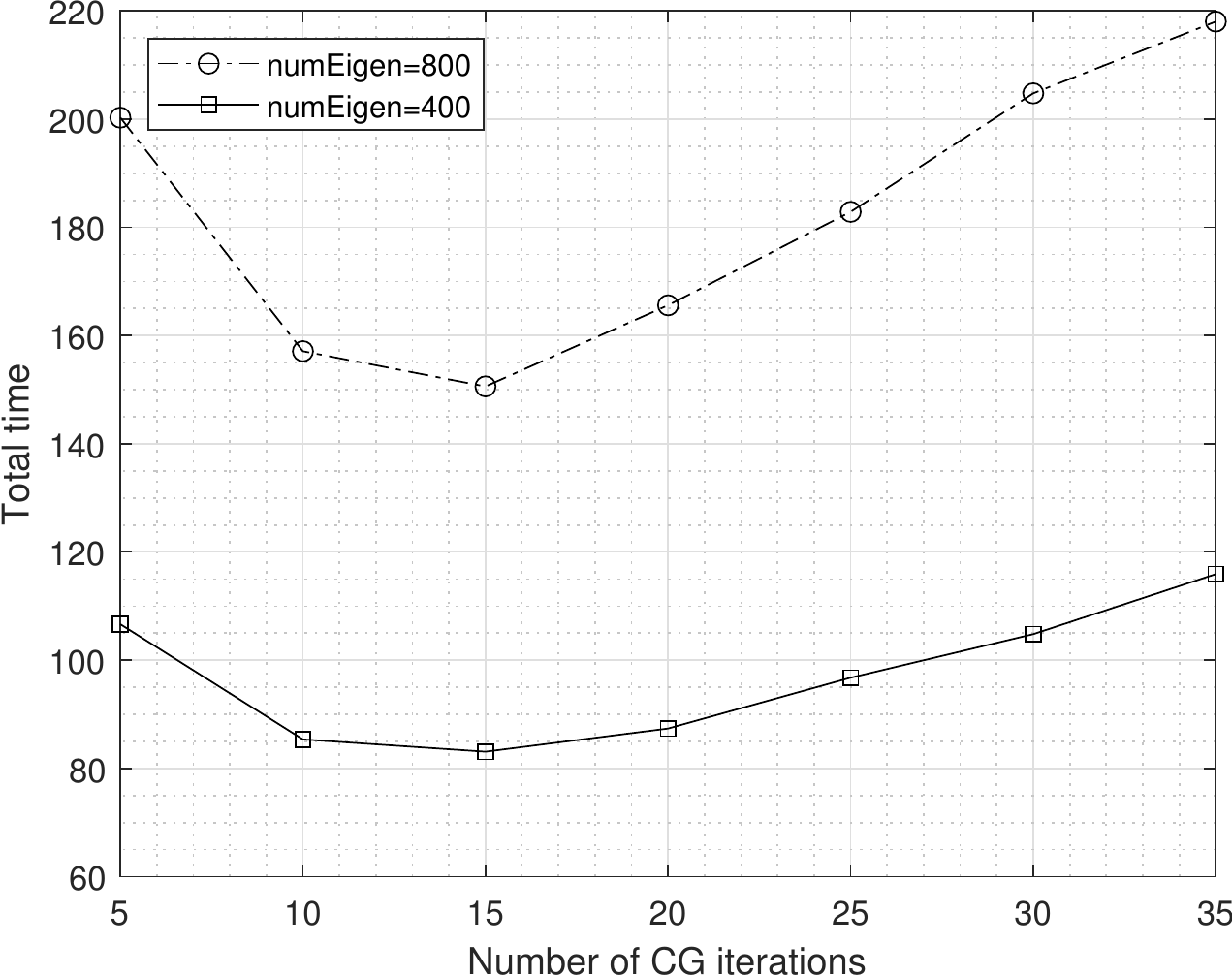}
\caption{Total time for SiO2 with ${\tt tol}=10^{-8}$}
\label{fig:cg_total_time_SiO2}
\end{figure}

For SiO2 with ${\tt numEigen}=400$ and $800$,
as shown in Figures \ref{fig:cg_nconv_SiO2} and \ref{fig:cg_avg_time_SiO2},
when the number of CG iterations is increased from $5$ to $35$ in each GCG iteration,
the number of GCG iterations decreases and the average time per each GCG iteration increases. 
And the total time reaches a minimum near $15$ CG iterations
according to Figure \ref{fig:cg_total_time_SiO2}.
In fact, from Andrews to SiO2, there have similar conclusions.

\begin{figure}[!htb]
\centering
\includegraphics[scale=0.45]{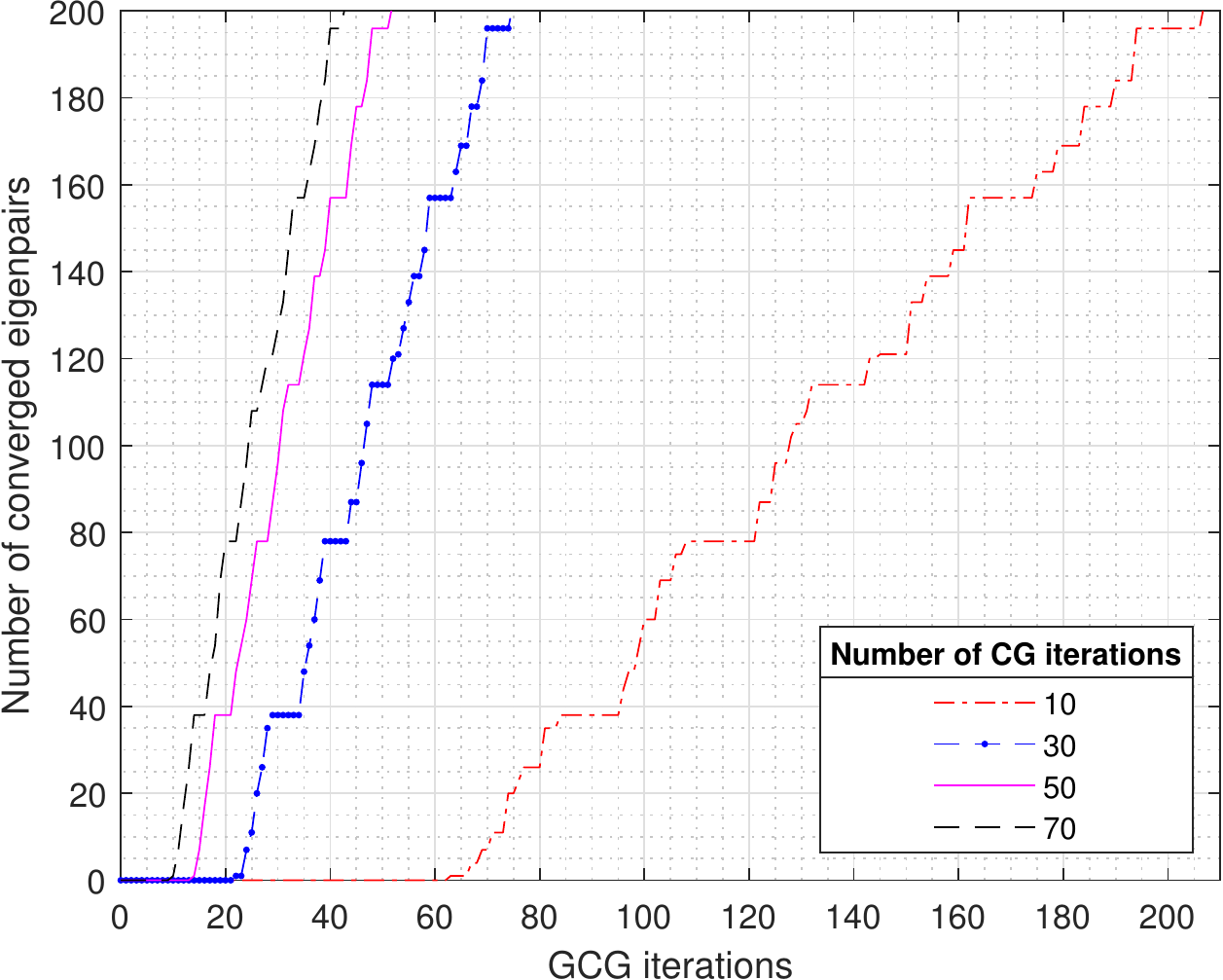}
\caption{Convergence procedure for FEM matrices with ${\tt tol}=10^{-8}$}
\label{fig:cg_nconv_200}
\end{figure}

\begin{figure}[!htb]
\centering
\includegraphics[scale=0.45]{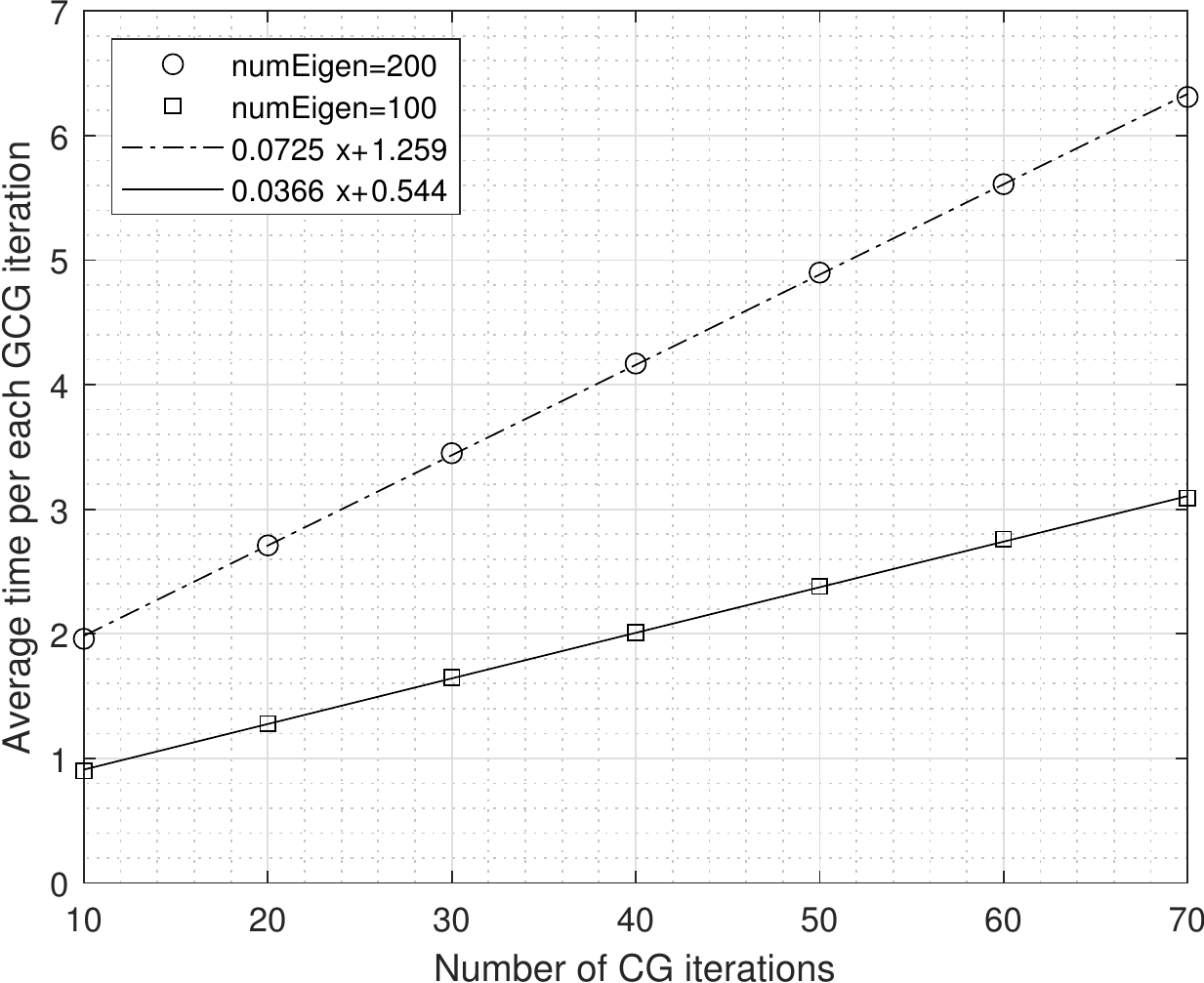}
\caption{Average time for FEM matrices with ${\tt tol}=10^{-8}$}
\label{fig:cg_avg_time}
\end{figure}

\begin{figure}[!htb]
\centering
\includegraphics[scale=0.45]{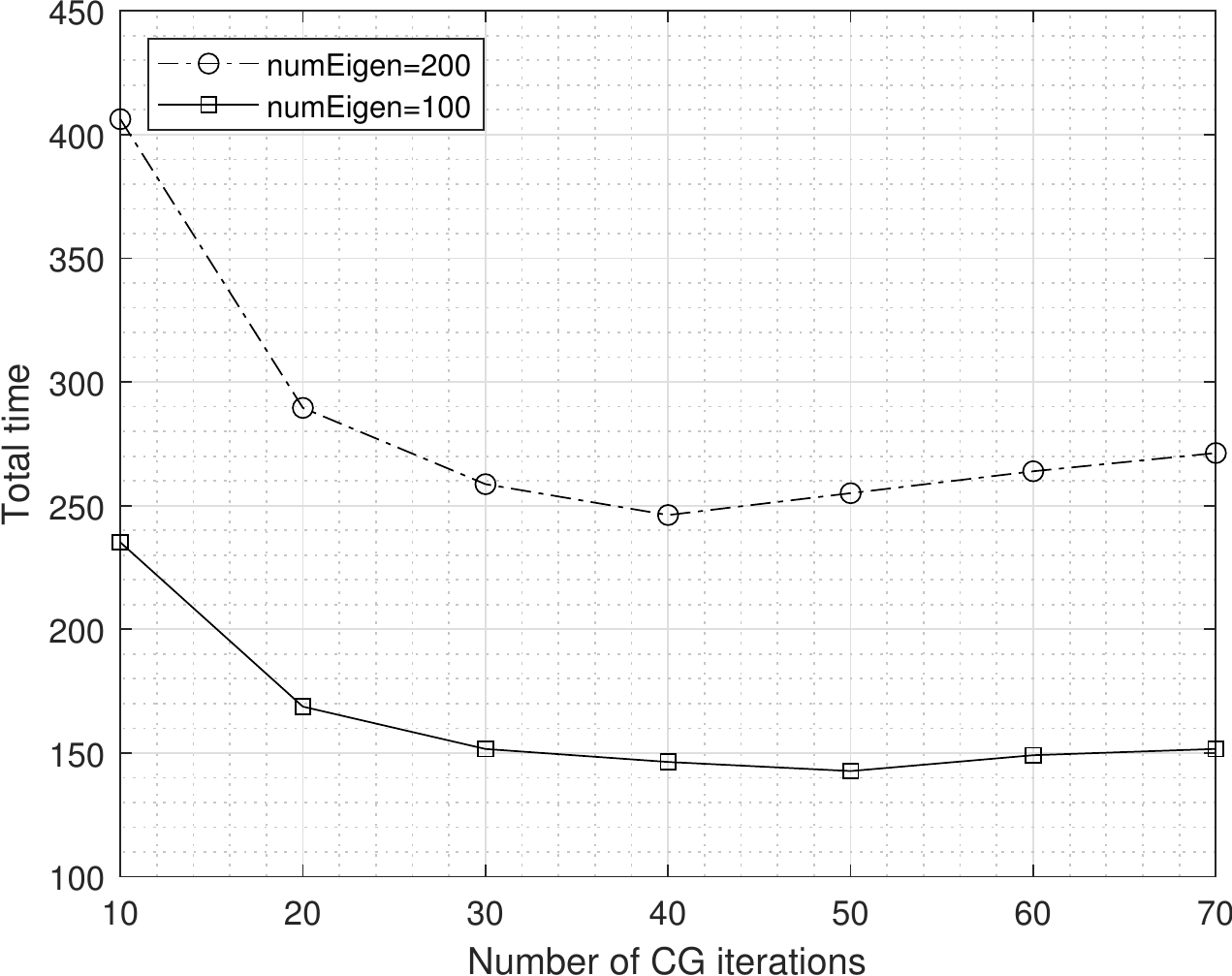}
\caption{Total time for FEM matrices with ${\tt tol}=10^{-8}$}
\label{fig:cg_total_time}
\end{figure}

Figures \ref{fig:cg_nconv_200} and \ref{fig:cg_avg_time} show the corresponding results for 
FEM matrices with ${\tt numEigen}=100$ and $200$.
When the number of CG iterations is increased from $10$ to $70$, 
the number of GCG iterations decreases and the average time per each GCG iteration increases.
The best performance is achieved at $30$-$40$ CG iterations
as shown in Figure \ref{fig:cg_total_time}.

It is noted that the number of CG iterations in each GCG iteration 
affects the efficiency of the algorithm deeply as presented
in Figures \ref{fig:cg_avg_time_SiO2} and \ref{fig:cg_avg_time}.
The average time per each GCG iteration 
is linearly associated with the number of CG iterations.
So, the number of CG iterations is a key parameter 
for trading off between the number of GCG iterations
and the average time of GCG iterations.
In fact, the total time of GCG algorithm is nearly equal to
the multiplication of
the number of GCG iterations
and
the average time of GCG iterations.
In other words,
though increasing the number of CG iterations
can accelerate convergence,
it takes more time in each GCG iteration.

In fact, the sparsity, the condition number 
and the dimension of the matrix
all affect the convergence rate of the CG iteration.
In the GCGE package, 
we set two stop conditions of the CG iteration.
When the residual of the solution is less than one percent of the initial residual, 
or the number of CG iterations is greater than $30$, 
the CG iteration will be stopped.
 

\subsection{About different tolerances}\label{sec:tol_small}

In this subsection, we will compare the performance of GCGE, LOBPCG and Krylov-Schur methods
under different tolerances.

\begin{figure*}[!htb]
\centering
\includegraphics[scale=0.45]{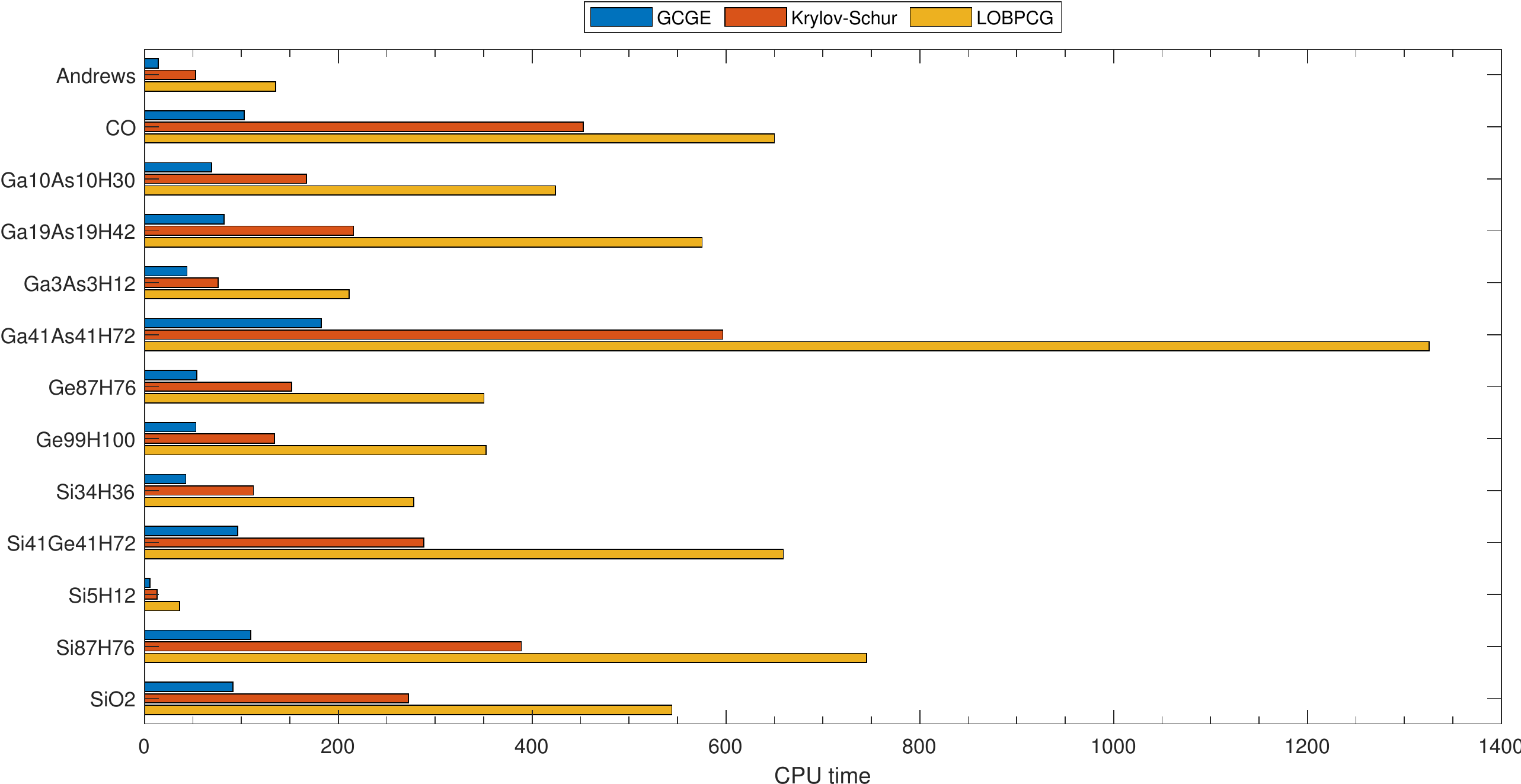}
\caption{${\tt tol}=10^{-4}$, ${\tt numEigen}=800$, and ${\tt numProc}=36$}
\label{fig:diff_method_1e-4}
\end{figure*}

\begin{figure*}[!htb]
\centering
\includegraphics[scale=0.45]{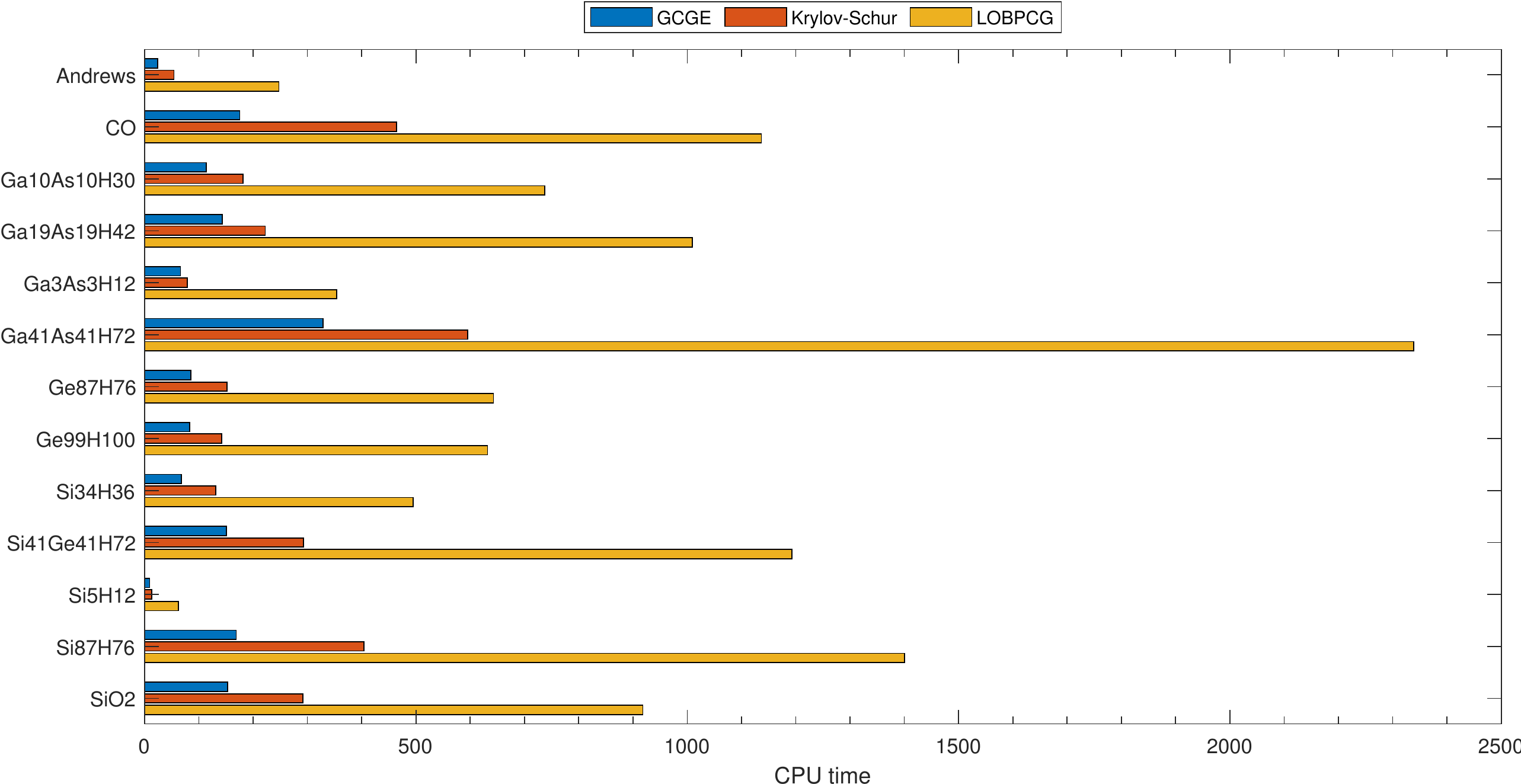}
\caption{${\tt tol}=10^{-8}$, ${\tt numEigen}=800$, and ${\tt numProc}=36$}
\label{fig:diff_method_1e-8}
\end{figure*}

\begin{figure*}[!htb]
\centering
\includegraphics[scale=0.45]{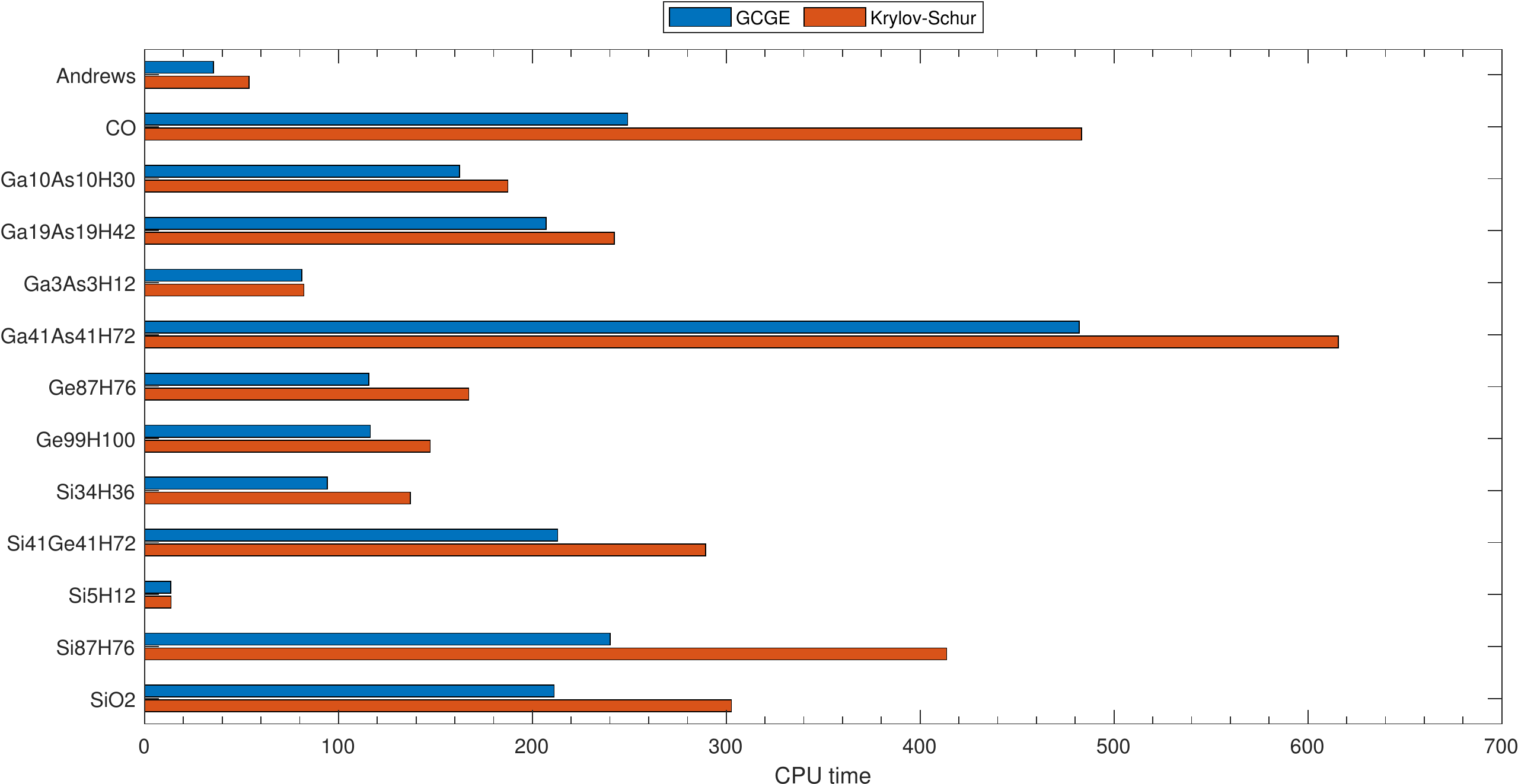}
\caption{${\tt tol}=10^{-12}$, ${\tt numEigen}=800$, and ${\tt numProc}=36$}
\label{fig:diff_method_1e-12}
\end{figure*}

In Figures \ref{fig:diff_method_1e-4} and \ref{fig:diff_method_1e-8}, 
GCGE, LOBPCG and Krylov-Schur methods with ${\tt numProc}=36$ are compared under 
${\tt tol}=10^{-4}$ and $10^{-8}$, respectively.
Under the tolerance $10^{-12}$,
LOBPCG can not converge after $3000$ iterations,
which means that the LOBPCG has no good stability.
So only the performance of GCGE and Krylov-Schur methods are 
compared under ${\tt tol}=10^{-12}$ and the results are 
presented in Figure \ref{fig:diff_method_1e-12}.
Here, MUMPS \cite{Amestoy2001fully,Amestoy2019Performance} is used as linear solver for Krylov-Schur method.

Obviously, GCGE is always more efficient than LOBPCG under different tolerances.
In addition,
when ${\tt tol}=10^{-4}$ and $10^{-8}$, 
GCGE is much faster than Krylov-Schur method.
Under tolerances $10^{-12}$, the CPU time of GCGE and Krylov-Shur method is similar and GCGE is slighly faster.

\begin{figure}[!htb]
\centering
\includegraphics[scale=0.45]{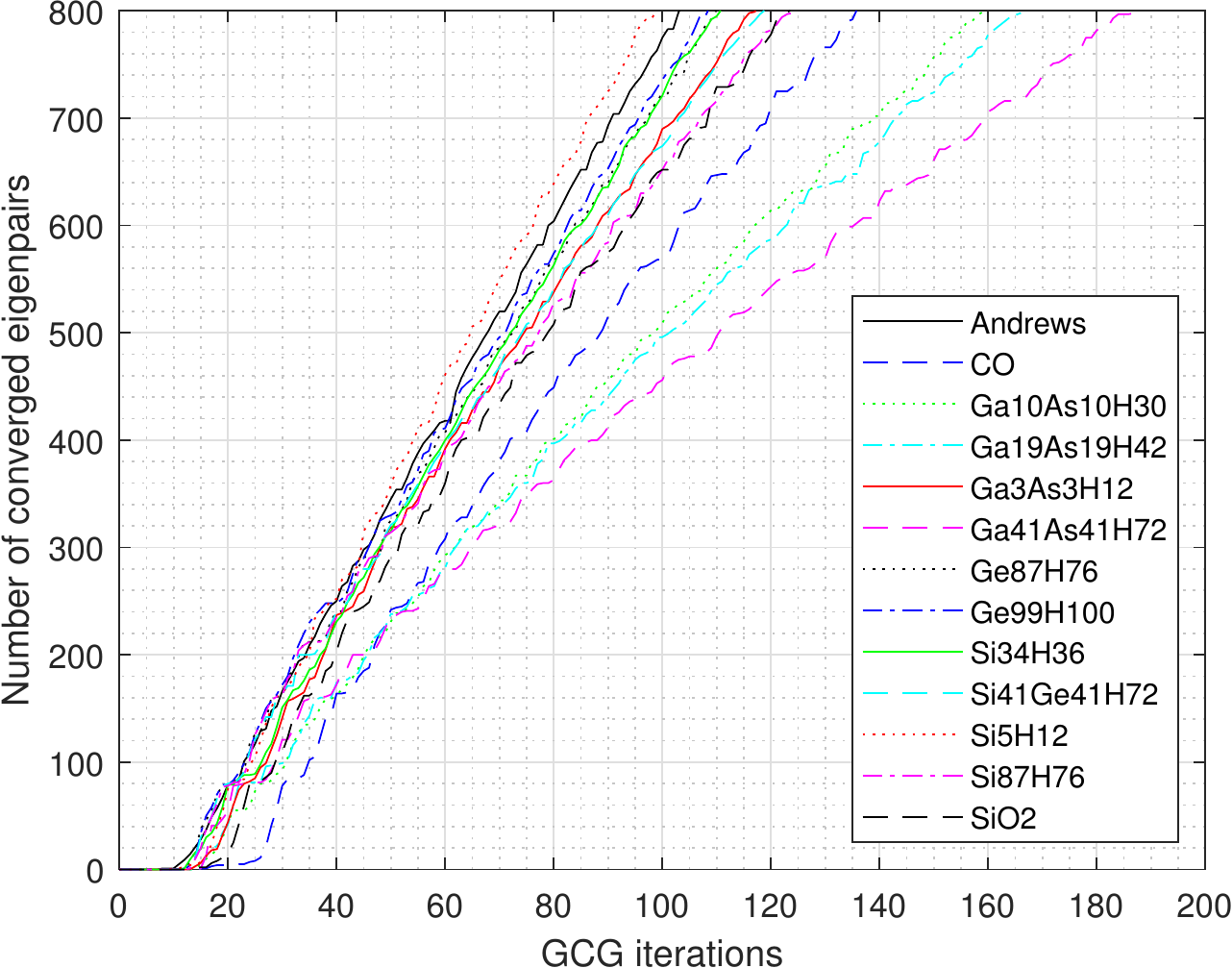}
\caption{${\tt tol}=10^{-12}$, ${\tt numEigen}=800$, and ${\tt numProc}=36$}
\label{fig:gcge_1e-12_ncon}
\end{figure}

\begin{figure}[!htb]
\centering
\includegraphics[scale=0.45]{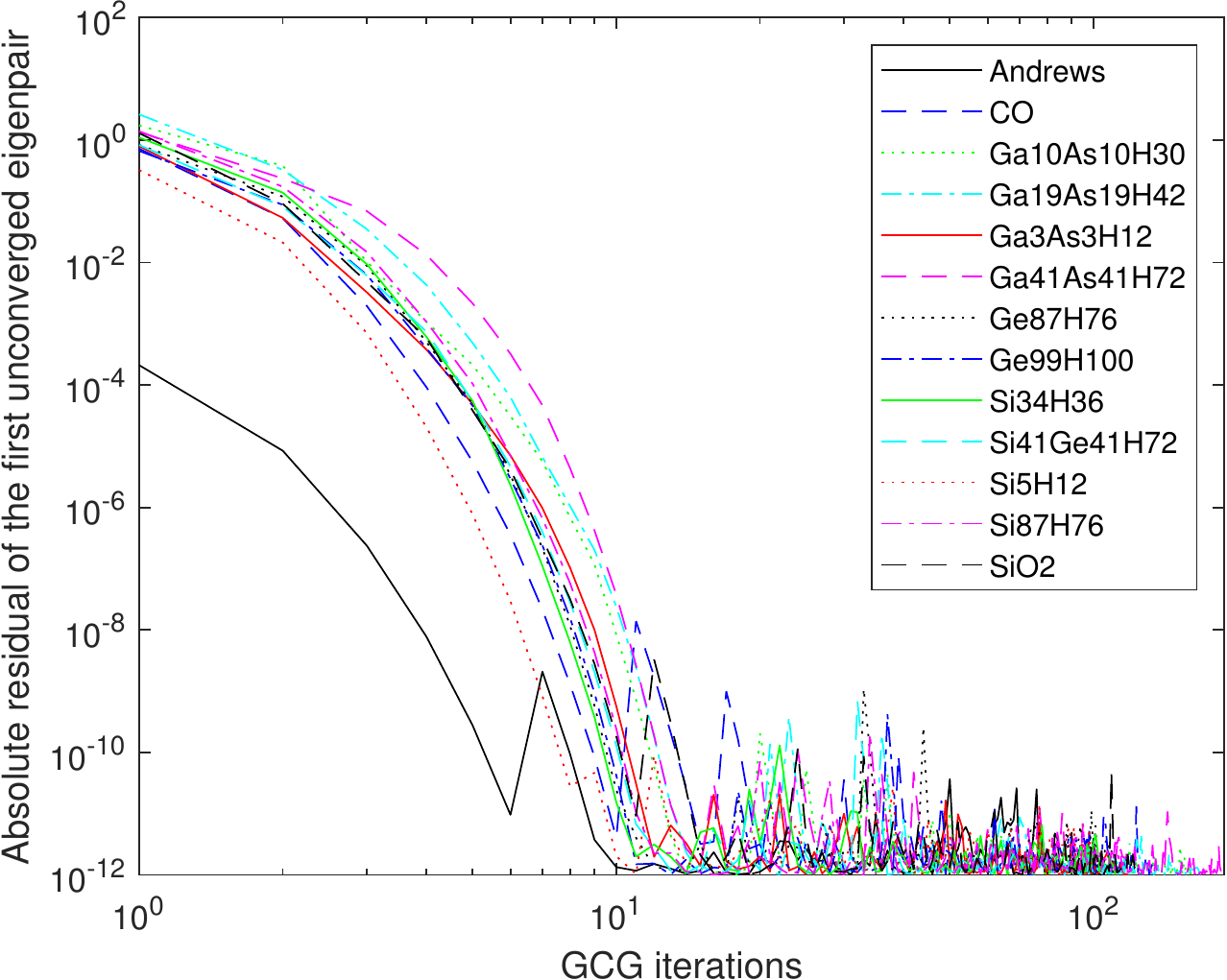}
\caption{${\tt tol}=10^{-12}$, ${\tt numEigen}=800$, and ${\tt numProc}=36$}
\label{fig:gcge_1e-12_resi}
\end{figure}

In addition,
the convergence procedure of GCG algorithm
with ${\tt tol}=10^{-12}$ 
for the first thirteen matrices
is shown 
in Figure \ref{fig:gcge_1e-12_ncon}.
As the number of GCG iterations increases,
the number of converged eigenpairs increases.
In Figure \ref{fig:gcge_1e-12_resi},
the absolute residual of the first unconverged eigenpair is presented.

\begin{figure}[!htb]
\centering
\includegraphics[scale=0.45]{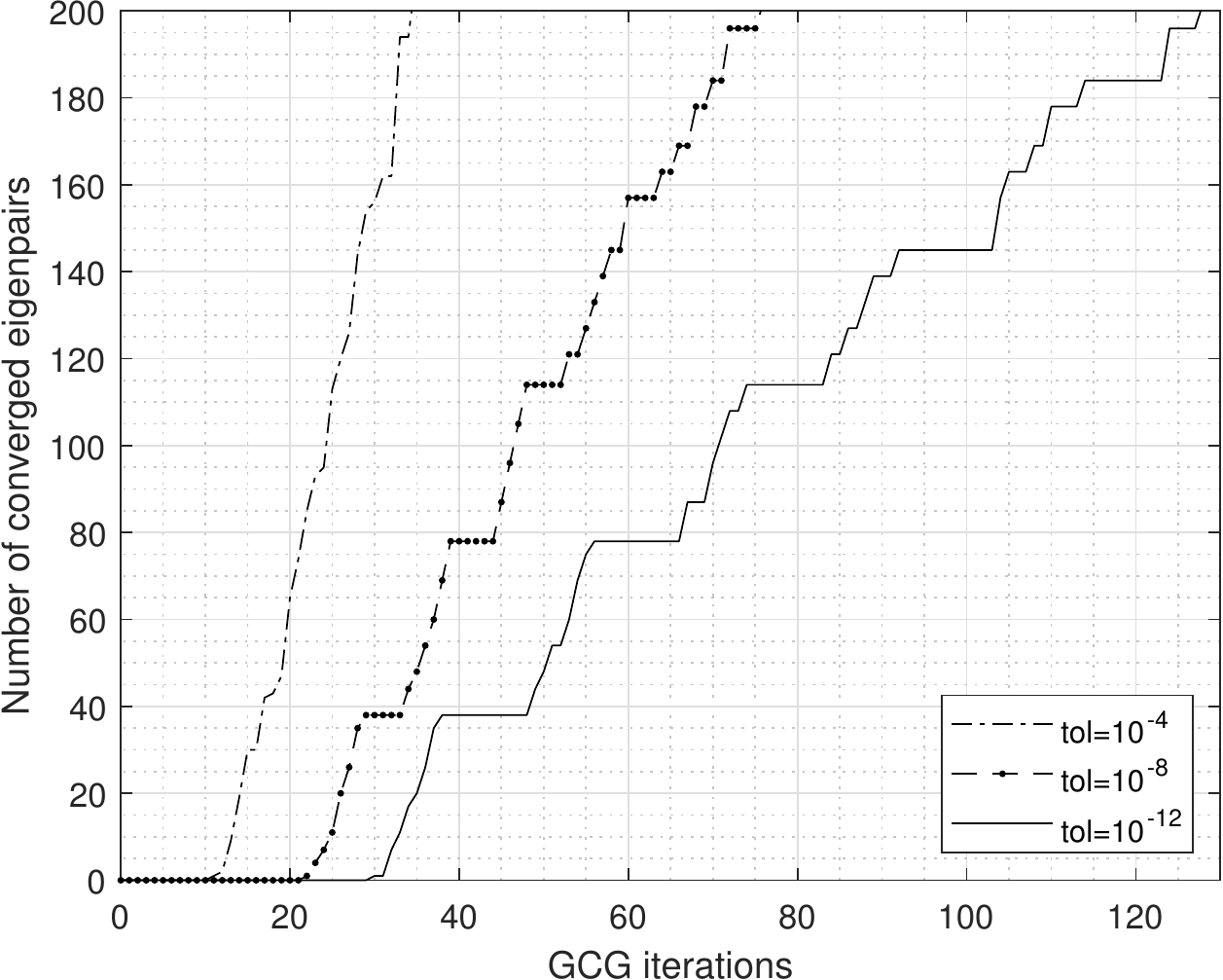}
\caption{${\tt numEigen}=200$ and ${\tt numProc}=576$}
\label{fig:gcge_ncon_FEM}
\end{figure}

\begin{figure}[!htb]
\centering
\includegraphics[scale=0.45]{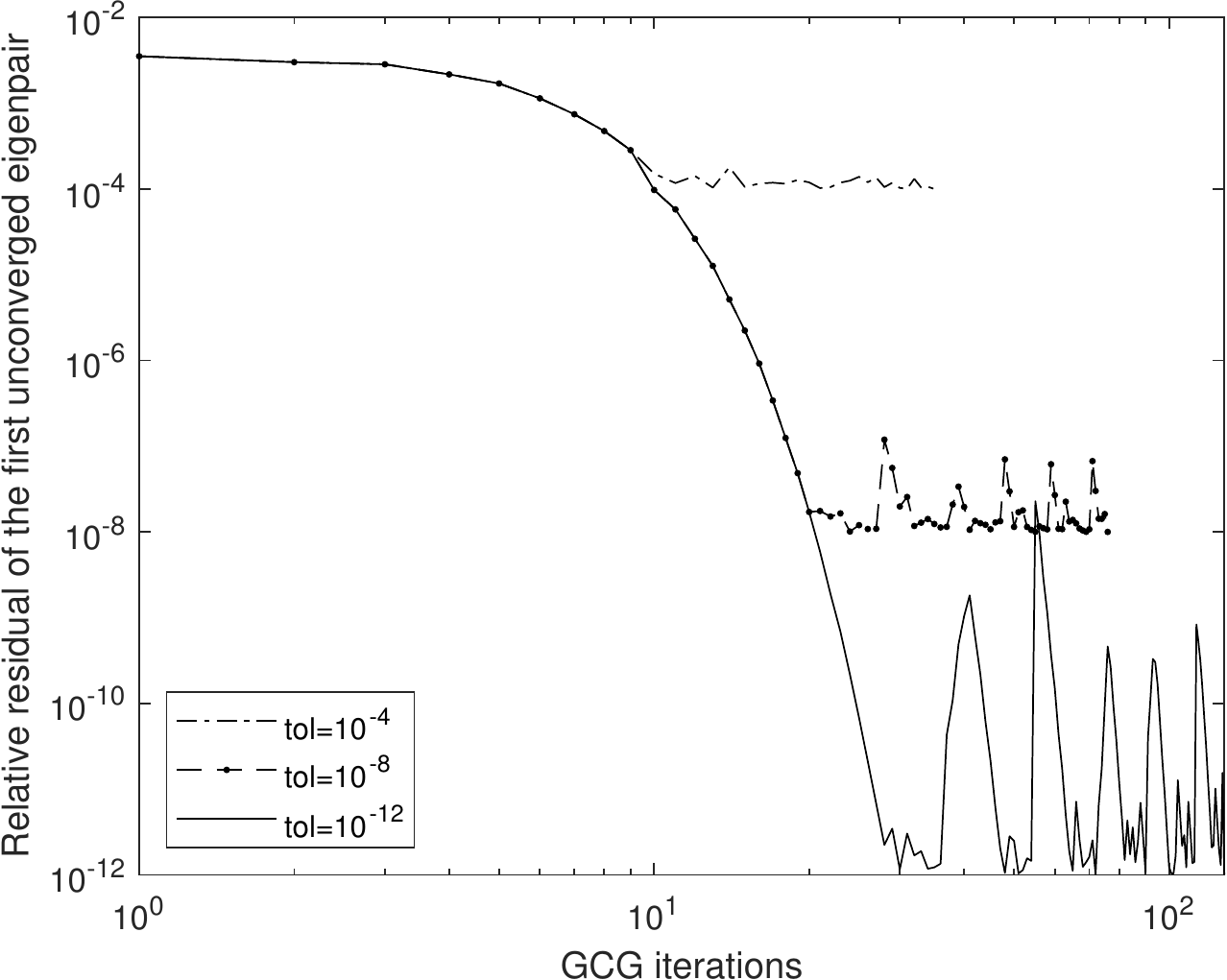}
\caption{${\tt numEigen}=200$ and ${\tt numProc}=576$}
\label{fig:gcge_resi_FEM}
\end{figure}

For FEM matrices, the performances of GCGE are shown 
in Figures \ref{fig:gcge_ncon_FEM}
and \ref{fig:gcge_resi_FEM}.
Due to ${\tt blockSize}={\tt numEigen}/5=40$,
there are four noticeable pauses for the case of ${\tt tol}=10^{-12}$
when the number of converged eigenpairs is close to 
$1\times40$, $2\times40$, $3\times40$, and $4\times40$
at around the 40th, 60th, 80th, and 100th GCG iteration.
Roughly speaking, 
the $40$ eigenpairs can be converged once every twenty GCG iterations.


\subsection{Scaling for the number of eigenpairs}

Here, we investigate the dependence of computing time
on the number of desired eigenpairs.
For this aim, we compute the first $50$-$800$ eigenpairs of matrices
listed in Table \ref{tab:matrix}.



\begin{figure}[!htb]
\centering
\includegraphics[scale=0.45]{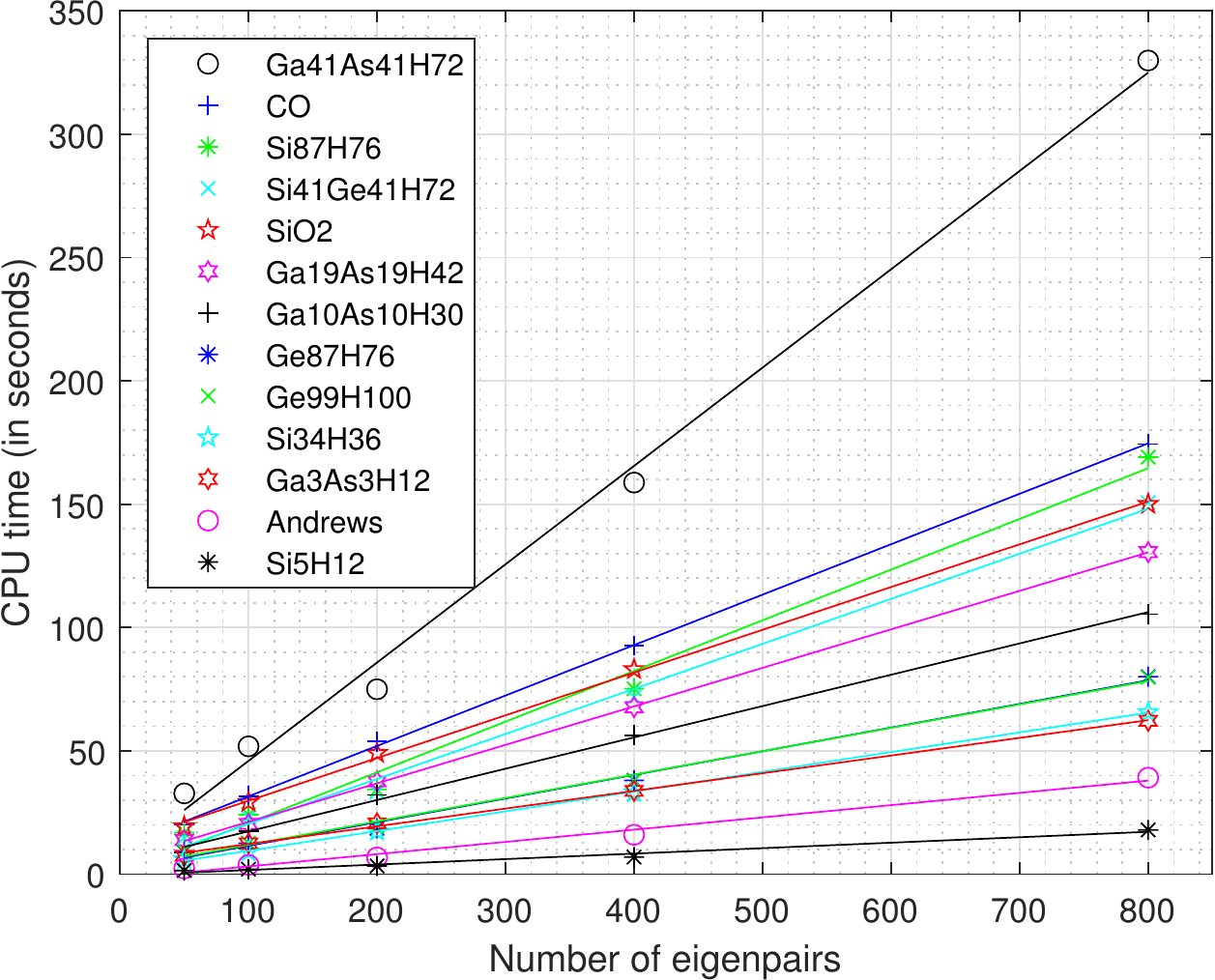}
\caption{GCGE with ${\tt tol}=10^{-8}$ and ${\tt numProc}=36$}
\label{fig:gcge_diff_nev}
\end{figure}

\begin{figure}[!htb]
\centering
\includegraphics[scale=0.45]{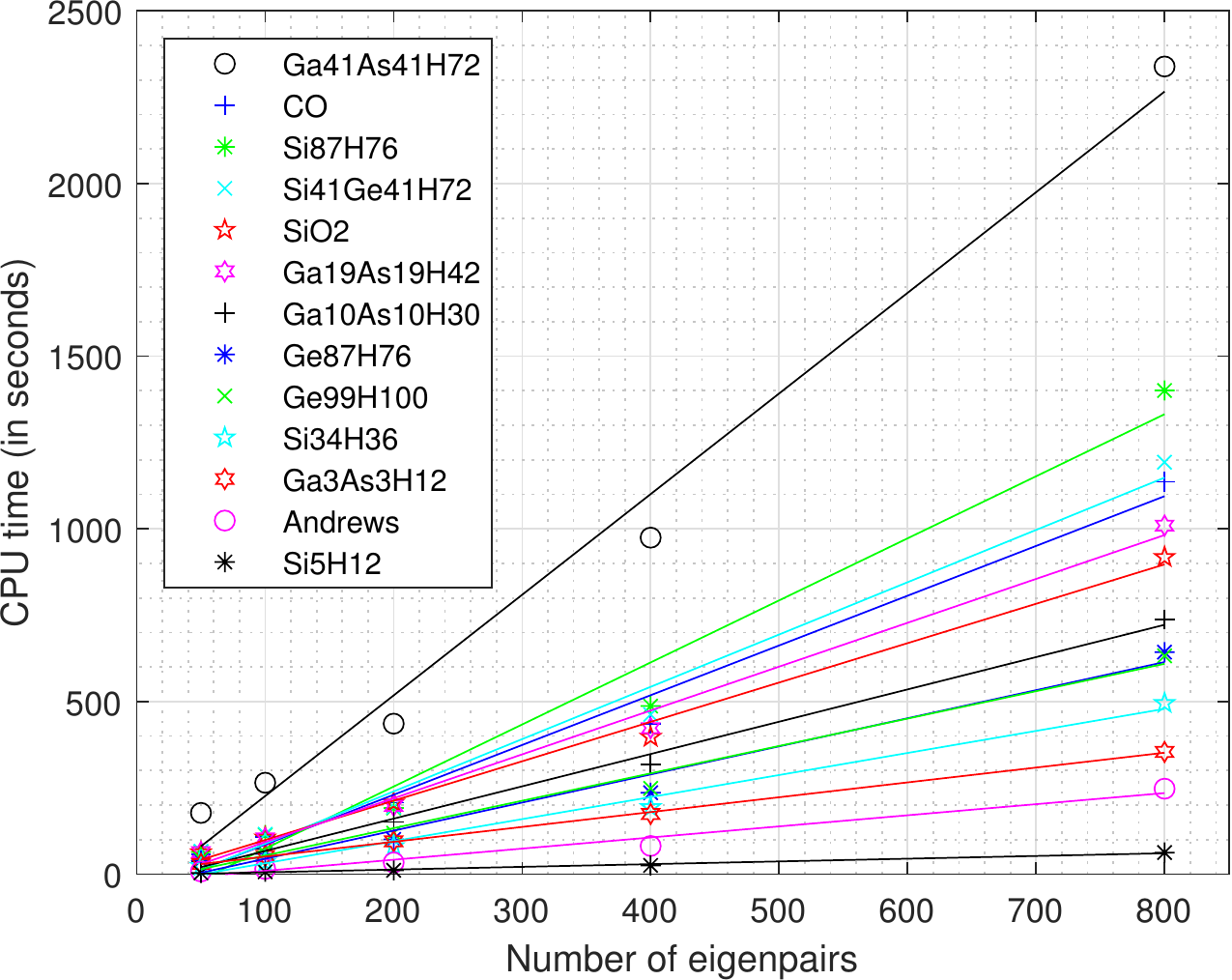}
\caption{LOBPCG with ${\tt tol}=10^{-8}$ and ${\tt numProc}=36$}
\label{fig:lobpcg_diff_nev}
\end{figure}

The test for the first thirteen matrices is performed on a single node with 36 processes.
The results in Figures \ref{fig:gcge_diff_nev} and \ref{fig:lobpcg_diff_nev}
show that
just like LOBPCG, GCGE has almost linear scaling property, 
which means the computing time is linearly dependent on the number of desired eigenpairs.
Moreover, GCGE has better efficiency than LOBPCG.
From Andrews to SiO2, the total time ratios of GCGE to LOBPCG are
\begin{align*}
&17.59\%,\ 19.17\%,\ 16.70\%,\ 15.02\%,\ 19.35\%,\ 15.46\%,\\    
&14.67\%,\ 14.43\%,\ 15.85\%,\ 14.44\%,\ 28.82\%,\ 14.15\%,\ 19.55\%.
\end{align*}

\begin{figure}[!htb]
\centering
\includegraphics[scale=0.45]{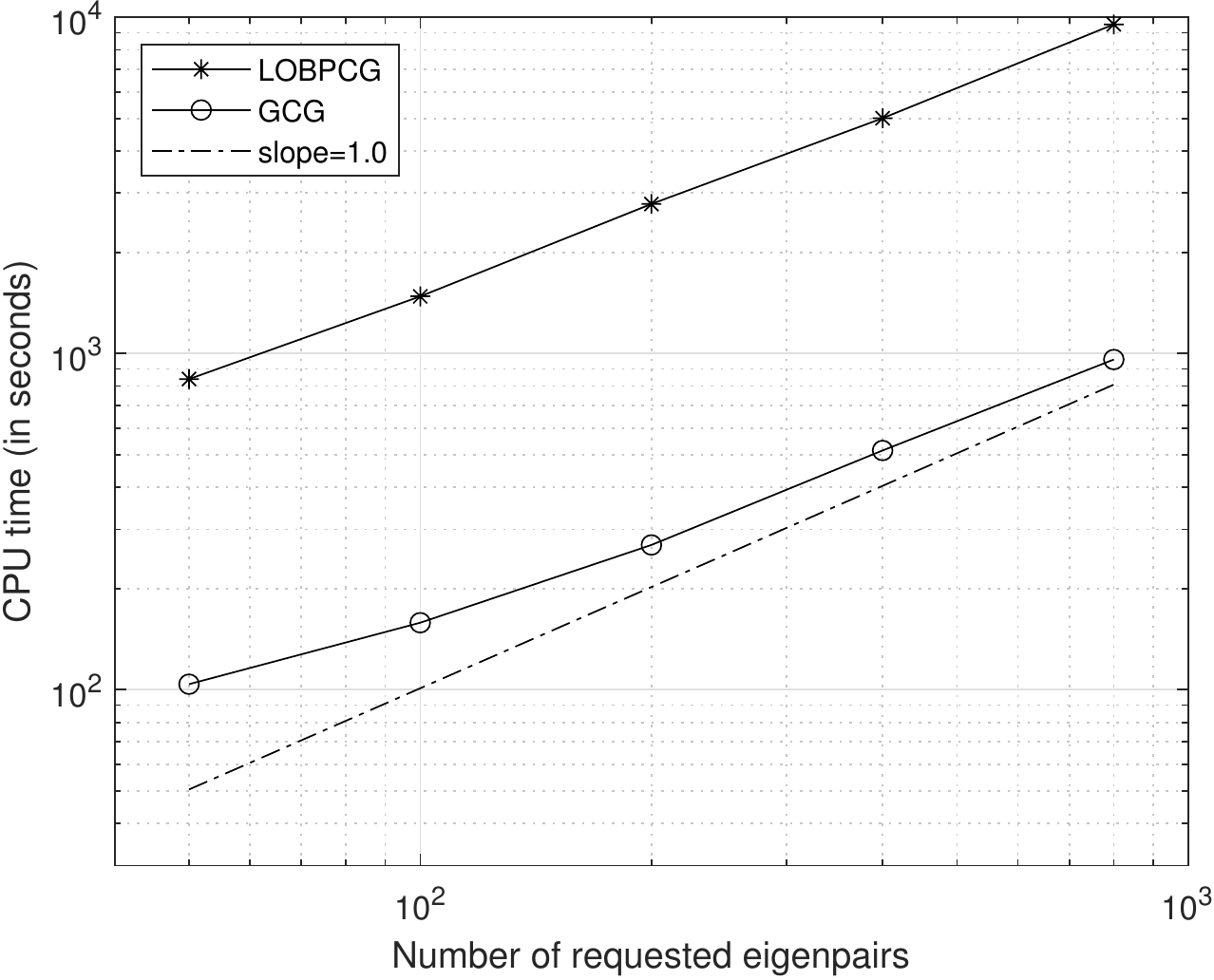}
\caption{CPU time for FEM matrices with ${\tt tol}=10^{-8}$ and ${\tt numProc}=576$}
\label{fig:cmp_nev}
\end{figure}

Since the scales of  FEM matrices are large, the test is performed with $576$ processes on $16$ nodes.
The dependence of CPU time (in seconds) for FEM matrices 
on the number of eigenpairs is shown in
Figure \ref{fig:cmp_nev},
which implies that GCGE has better efficiency than LOBPCG for large scale matrices.
Moreover, GCGE and LOBPCG both have almost linear scaling property for large scale FEM matrices.

\begin{remark}
In fact, Krylov-Schur method is low efficient for FEM matrices on multi-nodes.
In Table \ref{tab:gkl_fem_p1}, for ${\tt numEigen}=50,100,200$, 
the generalized eigenvalue problem is tested,
which is the discretization of the eigenvalue problem (\ref{Laplace_Eigenvalue_Problem})
for the conforming linear finite element (P1 element) 
with 3,145,728 elements.
The dimensions of the matrices $A$ and $B$ are both 512,191.
\begin{table}[!htb]
\centering
\caption{FEM matrices (P1 element) with ${\tt tol}=10^{-8}$ and ${\tt numProc}=36$}
\label{tab:gkl_fem_p1}
\begin{tabular}{rrrr}
	\toprule
Method      &   50    &   100   &    200\\
	\midrule                   
GCGE        &  20.15  &  38.98  &  71.49\\
Krylov-Schur&1032.33  &1360.56  &2180.28\\
LOBPCG      &  63.99  & 114.65  & 286.67\\
	\bottomrule
\end{tabular}
\end{table}
\end{remark}

\subsection{Scalability test}\label{subsec:scaling}

\begin{figure}[!htb]
\centering
\includegraphics[scale=0.45]{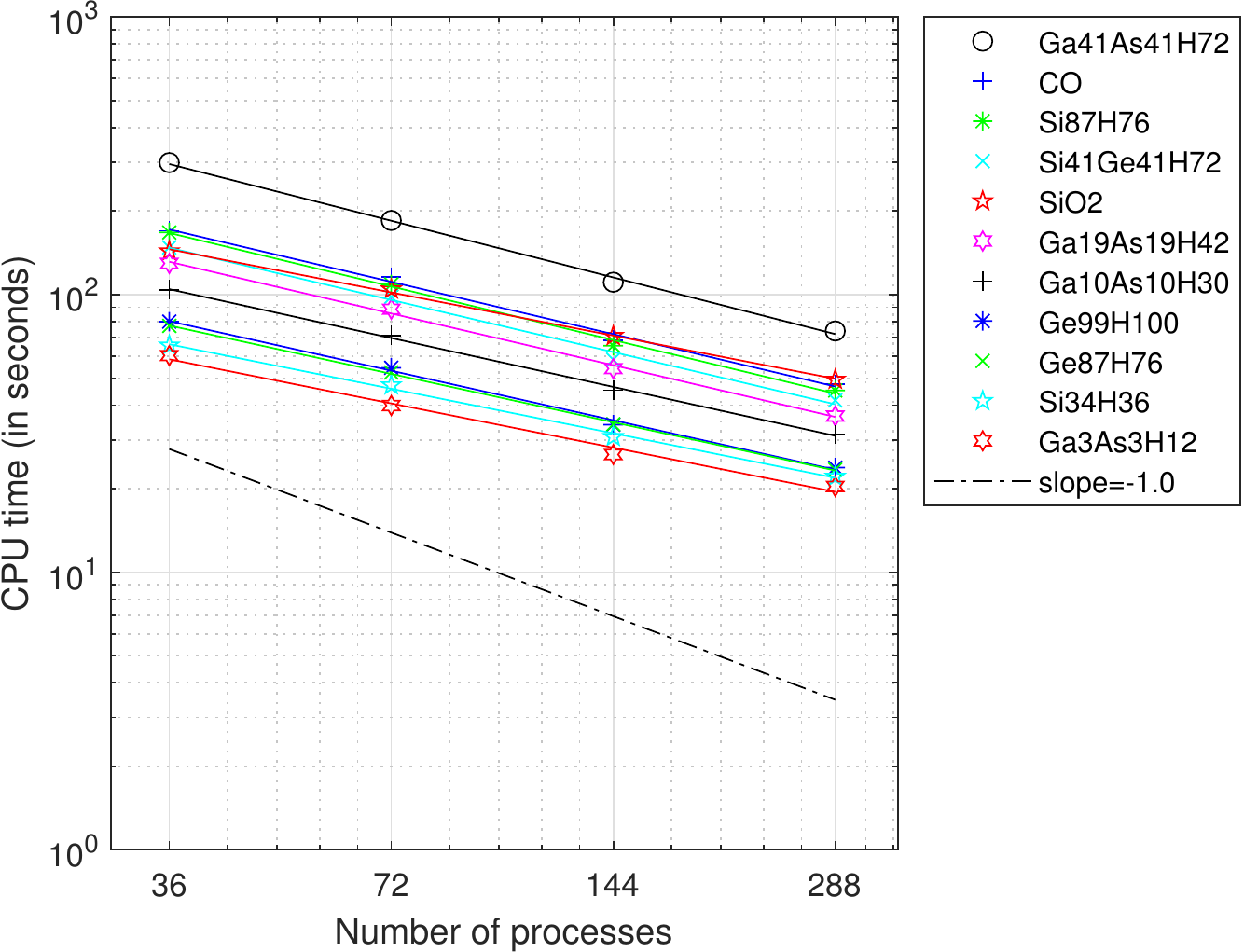}
\caption{GCGE with ${\tt tol}=10^{-8}$ and ${\tt numEigen}=800$}
\label{fig:gcge_diff_nproc}
\end{figure}

\begin{figure}[!htb]
\centering
\includegraphics[scale=0.45]{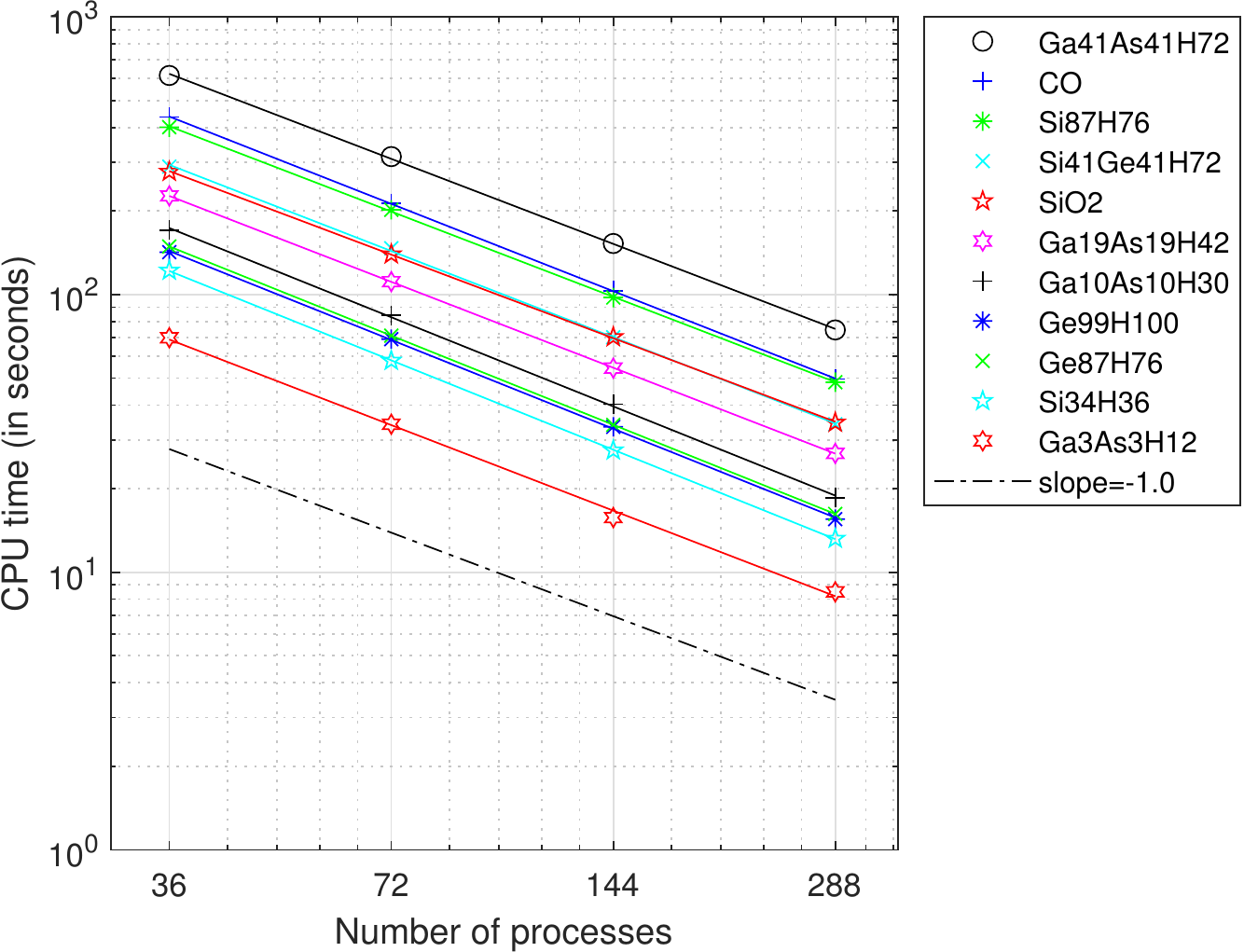}
\caption{Krylov-Schur method with ${\tt tol}=10^{-8}$ and ${\tt numEigen}=800$}
\label{fig:krylovschur_diff_nproc}
\end{figure}

\begin{figure}[!htb]
\centering
\includegraphics[scale=0.45]{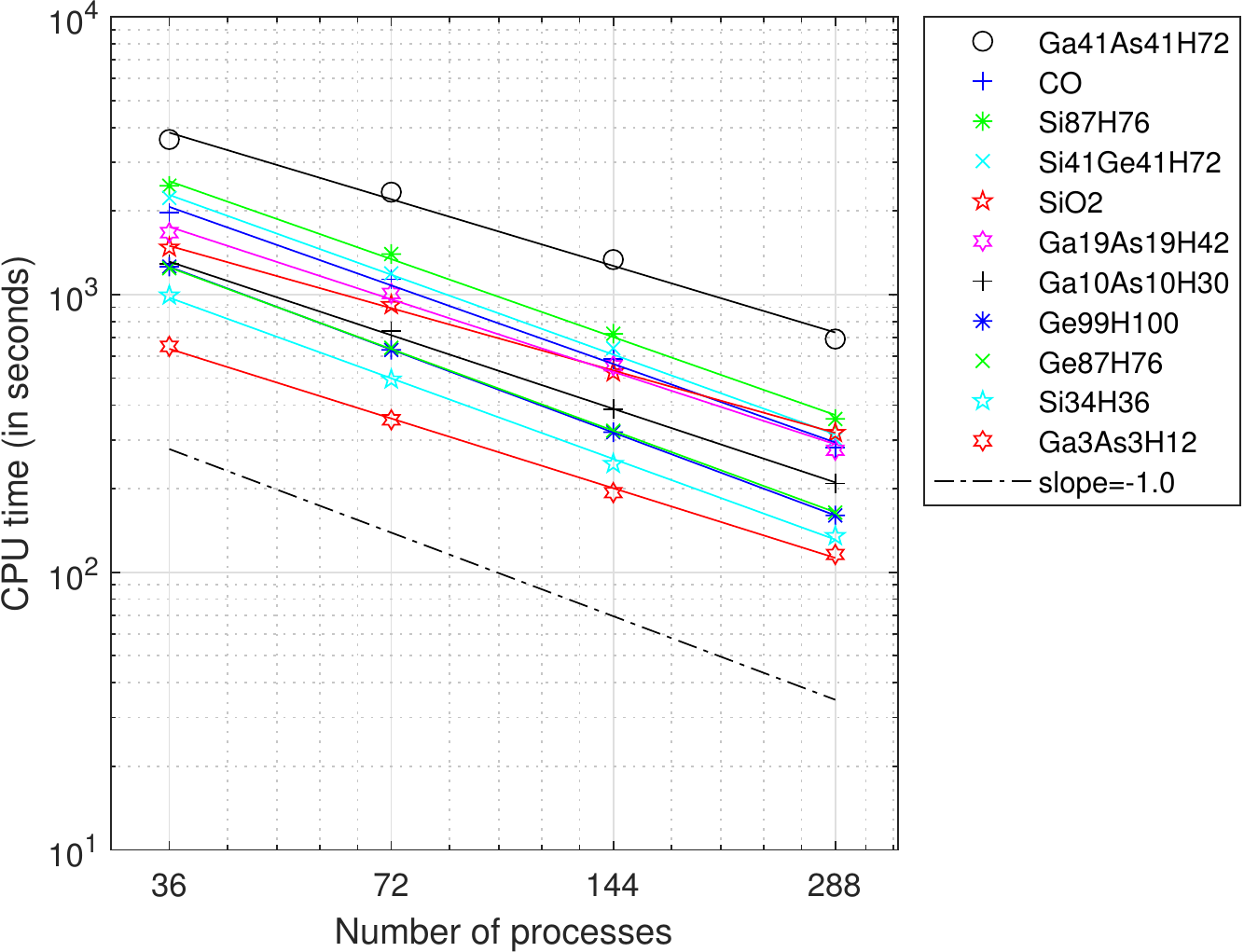}
\caption{LOBPCG with ${\tt tol}=10^{-8}$ and ${\tt numEigen}=800$}
\label{fig:lobpcg_diff_nproc}
\end{figure}

In order to do the scalability test,
we use $36$-$288$ processes to compute the first $800$ eigenpairs of the first thirteen matrices listed in Table \ref{tab:matrix}.
The comparisons of the scalability of GCGE and LOBPCG are shown
in Figures \ref{fig:gcge_diff_nproc}, \ref{fig:krylovschur_diff_nproc}, \ref{fig:lobpcg_diff_nproc}
and Table \ref{tab:gkl_a_g_s}.
It is noted that GCGE, LOBPCG, and Krylov-Schur methods have similar scalability 
for the first thirteen matrices, but 
the total time ratios of GCGE to LOBPCG are
\begin{align*}
&11.92\%,\ 10.10\%,\  9.61\%,\  8.79\%,\ 11.19\%,\  8.37\%,\\
& 7.88\%,\  8.10\%,\  8.86\%,\  7.93\%,\ 12.85\%,\  7.82\%,\ 11.36\%,
\end{align*}
from Andrews to SiO2.
In other words, GCGE has better efficiency than LOBPCG.
In addition, the total time ratios of GCGE to Krylov-Schur method are
\begin{align*}
&107.63\%,\ 50.08\%,\ 80.26\%,\ 73.64\%,\ 114.66\%,\ 57.77\%,\\
& 69.20\%,\ 73.86\%,\ 75.06\%,\ 64.07\%,\ 143.52\%,\ 51.67\%,\ 70.54\%,
\end{align*}
from Andrews to SiO2.
Only for small scale matrices Andrews (60,000), Ga3As3H12 (61,349), and Si34H36 (97,567),
the Krylov-Schur method is more efficient than GCGE, which are shown in Table \ref{tab:gkl_a_g_s}.

\begin{table}[!htb]
\centering
\caption{Small scale matrices with ${\tt tol}=10^{-8}$ and ${\tt numEigen}=800$}
\label{tab:gkl_a_g_s}
\begin{tabular}{rrrrr}
	\toprule
{\tt numProc} & Method &  Andrews &Ga3As3H12&   Si5H12\\
	\midrule                   
36 & GCGE        &  37.28 &  60.32&   9.38\\
   & Krylov-Schur&  54.66 &  69.90&  11.95\\
   & LOBPCG      & 447.09 & 650.75& 113.48\\
	\midrule               
72 & GCGE        &  26.09 &  39.85&   7.65\\
   & Krylov-Schur&  26.34 &  34.13&   5.30\\
   & LOBPCG      & 247.59 & 353.91&  62.59\\
	\midrule               
144& GCGE        &  22.69 &  26.54&   7.11\\
   & Krylov-Schur&  14.82 &  15.73&   2.90\\
   & LOBPCG      & 155.08 & 193.54&  44.10\\
	\midrule               
288& GCGE        &  34.55 &  20.36&   8.64\\
   & Krylov-Schur&  16.22 &   8.49&   2.69\\
   & LOBPCG      & 162.46 & 115.88&  34.98\\
	\bottomrule
\end{tabular}
\end{table}

\begin{figure}[!htb]
\centering
\includegraphics[scale=0.45]{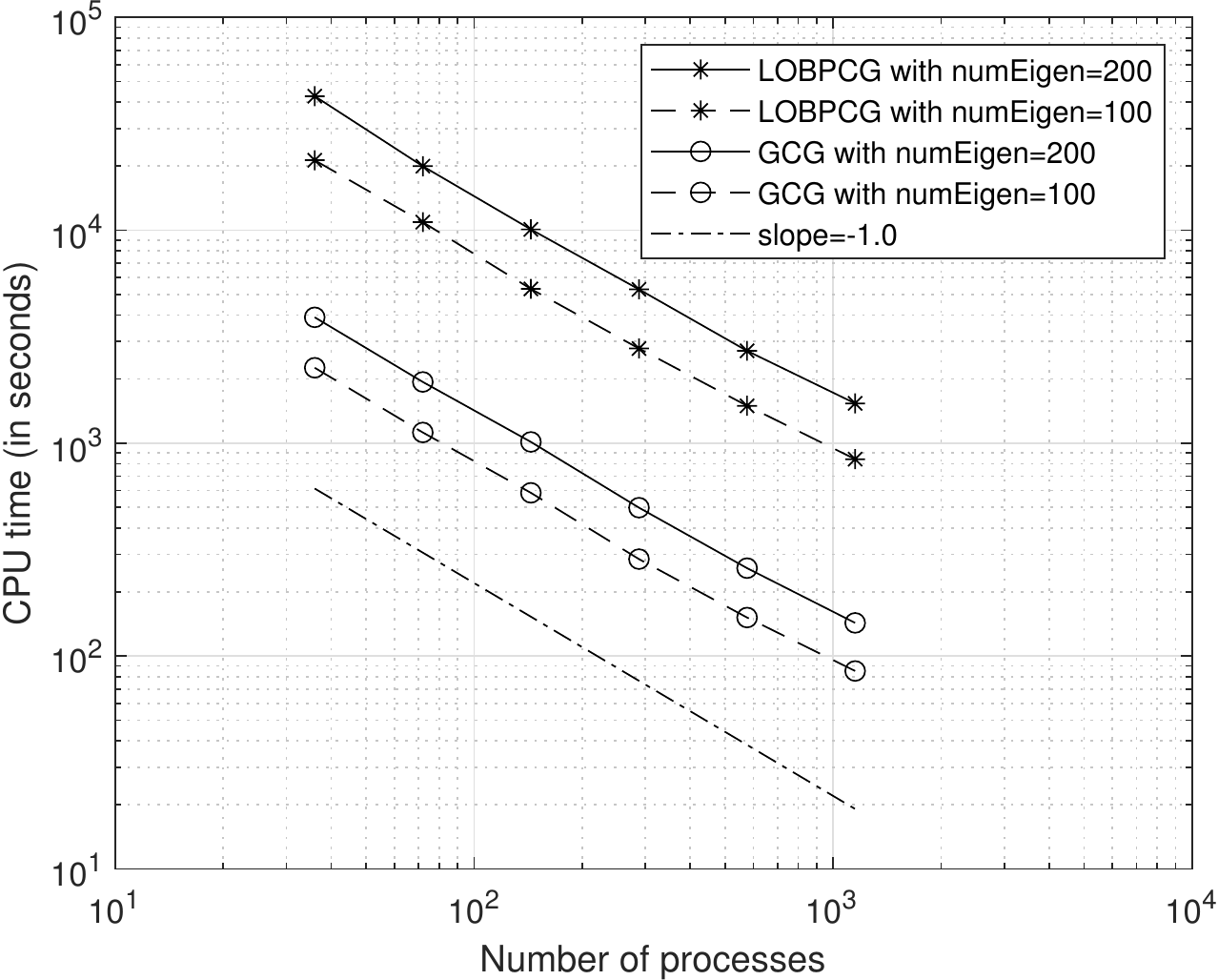}
\caption{CPU time for FEM matrices with ${\tt tol}=10^{-8}$}
\label{fig:cmp_fem}
\end{figure}

About the large scale FEM matrices, we use $36$-$1152$ processes 
for computing the lowest $100$ and $200$ eigenpairs.
In Figure \ref{fig:cmp_fem},
we can find that GCGE and LOBPCG have similar scalability
for large scale matrices, but
GCGE has better efficiency.
And the total time ratio of GCGE to LOBPCG is about $10\%$.

\subsection{The performance of GCGE with large {\tt numEigen}}\label{sec:nev_large}

In this subsection, 
the performance of the moving mechanism presented in Section \ref{sec:moving_mechanism} is tested.
The maximum project dimensions, ${\tt maxProjDim}$, are set to $1000$ and $2000$ for the first thirteen matrices and FEM matrices, respectively.

\begin{figure}[!htb]
\centering
\includegraphics[scale=0.45]{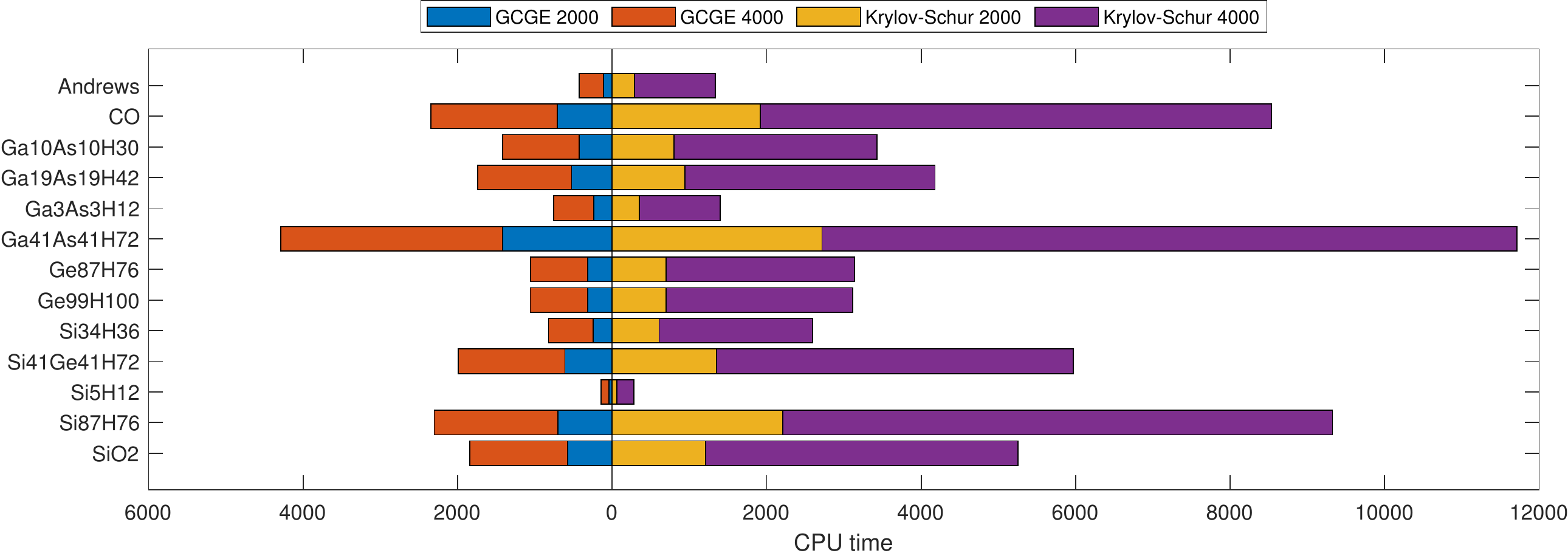}
\caption{${\tt tol}=10^{-12}$, ${\tt blockSize}=100$, and ${\tt numProc}=36$}
\Description{}
\label{fig:large_nev2000_4000}
\end{figure}

In Figure \ref{fig:large_nev2000_4000}, 
the performance of GCGE with the moving mechanism is shown 
for the first thirteen matrices,
For Krylov-Schur method,
we set ${\tt numEigen}$ to be $2000$ and $4000$
and the parameters are
\begin{verbatim}
    -eps_nev 2000 
    -eps_ncv 2400 
    -eps_mpd 800
\end{verbatim}
and 
\begin{verbatim}
    -eps_nev 4000 
    -eps_ncv 4400 
    -eps_mpd 1000
\end{verbatim}
respectively,
such that Krylov-Schur method has best efficiency for comparison.
Moreover, GCGE has better efficiency than Krylov-Schur.
From Andrews to SiO2, the total time ratios of GCGE to Krylov-Schur are 
\begin{align*}
&32.04\%,\ 27.49\%,\ 41.38\%,\ 41.70\%,\ 54.18\%,\ 36.60\%,\ 33.72\%,\\
&34.02\%,\ 31.76\%,\ 33.35\%,\ 50.08\%,\ 24.71\%,\ 35.05\%,\ 
\end{align*}

For FEM matrices with ${\tt numEigen}=5000$,
without the moving mechanism, 
the time of {\bf STEP 3} is dominated in Table \ref{tab:nev5000}.
And with the moving mechanism, 
the total time is reduced by about $50\%$.
In addition, 
the total time with dynamic shifts is reduced by about $20\%$ again
due to the reduction of the total number of GCG iterations.

\begin{table}[!htb]
\centering
\caption{The performance for FEM matrices
with ${\tt tol}=10^{-8}$, ${\tt blockSize}=200$, and ${\tt numProc}=1152$.}
\begin{tabular}{crr|rr|rr}
			\toprule
			&  \multicolumn{2}{c}{\mbox{Without Moving Mechanism}} & \multicolumn{2}{c}{\mbox{With Moving Mechanism}} & \multicolumn{2}{c}{\mbox{With Dynamic Shifts}}\\
		  \midrule
			&Time     & Percentage & Time & Percentage & Time & Percentage\\
		  \midrule
STEP 2& 445.05  & 5.33\%  & 42.64   & 1.04\%  & 41.56   & 1.25\% \\
STEP 3& 4727.57 & 56.65\% & 737.57  & 18.03\% & 601.92  &18.17\% \\
STEP 4& 78.94   & 0.95\%  & 34.09   & 0.83\%  & 32.60   & 0.98\% \\
STEP 5& 281.12  & 3.37\%  & 123.21  & 3.01\%  & 99.52   & 3.00\% \\
STEP 6& 2811.89 & 33.70\% & 3153.00 & 77.08\% & 2537.38 &76.59\% \\
		  \midrule              
			Total Time & 8344.57 & 100.00\% & 4090.52 &100.00\% & 3312.98 & 100.00\% \\
		  \midrule              
			Ratio      & 100.00\%&          & 49.02\% &         & 39.70\% &          \\
			\bottomrule
\end{tabular}
	\label{tab:nev5000}
\end{table}

In Table \ref{tab:nev10000_FEM},
the performances of two different orthogonalization methods are also compared.
When {\tt numEigen = 10000},
Algorithm \ref{RecusiveOrthSVD} is faster than 
Algorithm \ref{Modified_Block_Orth}
because of fewer multiplication of matrix and vectors,
especially for the generalized algebraic eigenvalue problems.

\begin{table}[!htb]
\centering
\caption{The performance for FEM matrices
with ${\tt tol}=10^{-8}$, ${\tt blockSize}=200$, and ${\tt numProc}=1152$.}
\begin{tabular}{crr|rr}
			\toprule
			&  \multicolumn{2}{c}{Algorithm \ref{Modified_Block_Orth}} &\multicolumn{2}{c}{Algorithm \ref{RecusiveOrthSVD}} \\
		  \midrule
			&Time     & Percentage & Time & Percentage \\
		  \midrule
STEP 2& 18.37    & 0.16\%  & 18.35 	 & 0.17\%  \\
STEP 3& 1992.81  & 17.58\% & 2025.87 & 18.45\% \\
STEP 4& 90.04    & 0.79\%  & 66.93 	 & 0.61\%  \\
STEP 5& 324.76   & 2.87\%  & 326.45  & 2.97\%  \\
STEP 6& 8907.10  & 78.59\% & 8544.80 & 77.80\% \\
		  \midrule
Total Time & 11333.08 & 100.00\% & 10982.40 & 100.00\% \\
		  \midrule
			Ratio      & 100.00\% &          & 96.90\% &          \\
			\bottomrule
\end{tabular}
	\label{tab:nev10000_FEM}
\end{table}

\section{Concluding remarks}
This paper highlights some new issues for computing plenty of eigenpairs of large scale matrices on high performance computers.
The GCGE package is presented which is built 
with the damping block inverse power method with dynamic shifts for symmetric eigenvalue problems. 
Furthermore, in order to improve the efficiency, stability and scalability 
of the concerned package, the new efficient implementing techniques are designed
for updating subspaces, orthogonalization and computing Rayleigh-Ritz problems. 
Plenty of numerical tests are provided to validate the proposed package GCGE,
which can be downloaded from
\url{https://github.com/Materials-Of-Numerical-Algebra/GCGE}.



\begin{acks}
This research is supported partly by 
National Key R\&D Program of China 2019YFA0709600, 2019YFA0709601,
National Natural Science Foundations of China (Grant No. 11771434),
the National Center for Mathematics and Interdisciplinary Science, CAS,
and
Tianjin Education Commission Scientific Research Plan (2017KJ236).
\end{acks}

\bibliographystyle{ACM-Reference-Format}
\bibliography{F:/liyu/fullbib/fullbib}


\end{document}
\endinput